\documentclass[12pt]{article}
\usepackage[left=2cm,top=2cm,right=2cm,bottom=2cm]{geometry}
\usepackage[british]{babel}
\usepackage{amsfonts,amssymb,amsmath,latexsym,hyperref,fancyhdr}
\usepackage{graphicx}
\usepackage{epstopdf}
\usepackage{algorithmic}
\usepackage[british]{babel}
\usepackage{enumerate,bbm,xcolor,theorem,enumitem,verbatim,stmaryrd}

\pdfoutput=1

\ifpdf
  \DeclareGraphicsExtensions{.eps,.pdf,.png,.jpg}
\else
  \DeclareGraphicsExtensions{.eps}
\fi

\catcode`\<=\active \def<{
\fontencoding{T1}\selectfont\symbol{60}\fontencoding{\encodingdefault}}
\catcode`\>=\active \def>{
\fontencoding{T1}\selectfont\symbol{62}\fontencoding{\encodingdefault}}
\newcommand{\mathd}{\mathrm{d}}
\newcommand{\nin}{\not\in}

\newcommand{\tmem}[1]{{\em #1\/}}
\newcommand{\tmmathbf}[1]{\ensuremath{\boldsymbol{#1}}}
\newcommand{\tmop}[1]{\ensuremath{\operatorname{#1}}}

{\theorembodyfont{\rmfamily}\newtheorem{remark}{Remark}}
{\theorembodyfont{\rmfamily}\newtheorem{lemma}{Lemma}}

\newcommand\smallO{
  \mathchoice
    {{\scriptstyle\mathcal{O}}}
    {{\scriptstyle\mathcal{O}}}
    {{\scriptscriptstyle\mathcal{O}}}
    {\scalebox{.6}{$\scriptscriptstyle\mathcal{O}$}}
  }

\numberwithin{equation}{section}

\definecolor{myred}{RGB}{160,0,0}
\definecolor{mygreen}{RGB}{0,160,0}
\definecolor{myblue}{RGB}{0,0,160}
\hypersetup{
    colorlinks = true,
    linkcolor = {myred},
    citecolor = {mygreen},
    urlcolor  = {myblue},
}

\newcommand{\RED}{} 

\definecolor{manchester}{rgb}{.42,.17,.58}
\newcommand{\PUR}{} 

\title{On the asymptotic properties\\ of \RED{a canonical} diffraction integral}

\author{Rapha\"{e}l C. Assier$^{*}$ and I. David Abrahams$^{\dagger}$\\
\footnotesize{$^{*}$ Department of Mathematics, University of Manchester, Oxford Road, Manchester, {\rm M13 9PL}, UK}\\
\footnotesize{$^{\dagger}$ Isaac Newton Institute, University of Cambridge, 20 Clarkson Road, Cambridge CB3 0EH, UK}
}

\begin{document}

\maketitle

\begin{abstract}
We introduce and study a new \RED{canonical} integral, denoted $I_{+
  -}^{\varepsilon}$, depending on two complex parameters $\alpha_1$ and
  $\alpha_2$. It arises from the problem of wave diffraction by a
  quarter-plane, and is heuristically constructed to capture the complex field near
  the tip and edges. We establish some region of analyticity of this integral
  in $\mathbb{C}^2$, and derive its rich asymptotic behaviour as $| \alpha_1
  |$ and $| \alpha_2 |$ tend to infinity. We also study the decay properties
  of the function obtained from applying a specific double Cauchy integral
  operator to this integral. These results allow us to show that this integral
  shares all of the asymptotic properties expected from the key unknown function $G_{+
  -}$ arising when the quarter-plane diffraction problem is studied via a
  two-complex-variables Wiener--Hopf technique (see Assier \& Abrahams,
  \href{https://arxiv.org/abs/1905.03863}{arXiv:1905.03863}, 2020). As a result, the integral $I_{+ -}^{\varepsilon}$
  can be used to mimic the unknown function $G_{+ -}$ and to build an
  efficient `educated' approximation to the quarter-plane problem. 

\end{abstract}

\section{Introduction and motivation}
We propose to study the properties (asymptotic behaviour, analyticity and more) of the integral $I_{+-}^\varepsilon$ that can be considered a function of two complex variables $(\alpha_1,\alpha_2)$ and is defined by
\begin{align}
I_{+ -}^{\varepsilon} (\alpha_1, \alpha_2) =  \int_{\frac{3 \pi}{2}}^{2
  \pi} \frac{\RED{f}\left( \frac{\pi}{2}, \varphi \right)}{(- i (\alpha_1 \cos
  (\varphi) + \alpha_2 \sin (\varphi) + i \varepsilon))^{\RED{\nu} + 3 / 2}}
  \mathd \varphi , \label{eq:initialdefIrevision}
\end{align}
where the constants $\varepsilon$ and $\RED{\nu}$ and the function $\RED{f}\left(\frac{\pi}{2}, \varphi \right)$ will be specified later. 
This work is directly motivated by the conclusions discussed in a recent article by the authors, {\cite{Assier2019}}, which focuses on a two-complex-variable investigation of the
three-dimensional problem of wave diffraction by a quarter-plane with homogeneous Dirichlet boundary conditions subject to an
incident plane wave. 

\RED{The quarter-plane problem is an unsolved \textit{Canonical Scattering Problem} (CSP). Typically, CSPs relate to wave diffraction by `simple' geometries exhibiting some characteristic features such
as curvature, edges or corners. They are the building blocks of most wave diffraction approximation techniques (such as the Geometrical Theory of Diffraction \cite{Keller1962}) used for more complicated, real-life, geometries and are extremely important in applications including e.g. noise
reduction and radar. Famous examples of solved CSPs are e.g. the Sommerfeld problem of diffraction by a half-plane \cite{Sommerfeld1896} and wedge diffraction problems; both of these two-dimensional CSPs can be tackled by the Wiener-Hopf technique in one complex variable 
\cite{Noble1958,Lawrie2007,Nethercote2020}. The solutions are expressed as complex integrals whose far-field asymptotic behaviour can be written down exactly. The quarter-plane problem can be seen as the direct three-dimensional generalisation of the Sommerfeld problem.}

\RED{Another well-known canonical integral, akin to (\ref{eq:initialdefIrevision}) in the sense that it can be interpreted as a function of two complex variables, and with application to scattering amplitude in string theory \cite{Veneziano1968}, is the Beta function (also known as the Euler integral of the first kind).}

\begin{figure}[h]
\centering
  \includegraphics[width=0.48\textwidth]{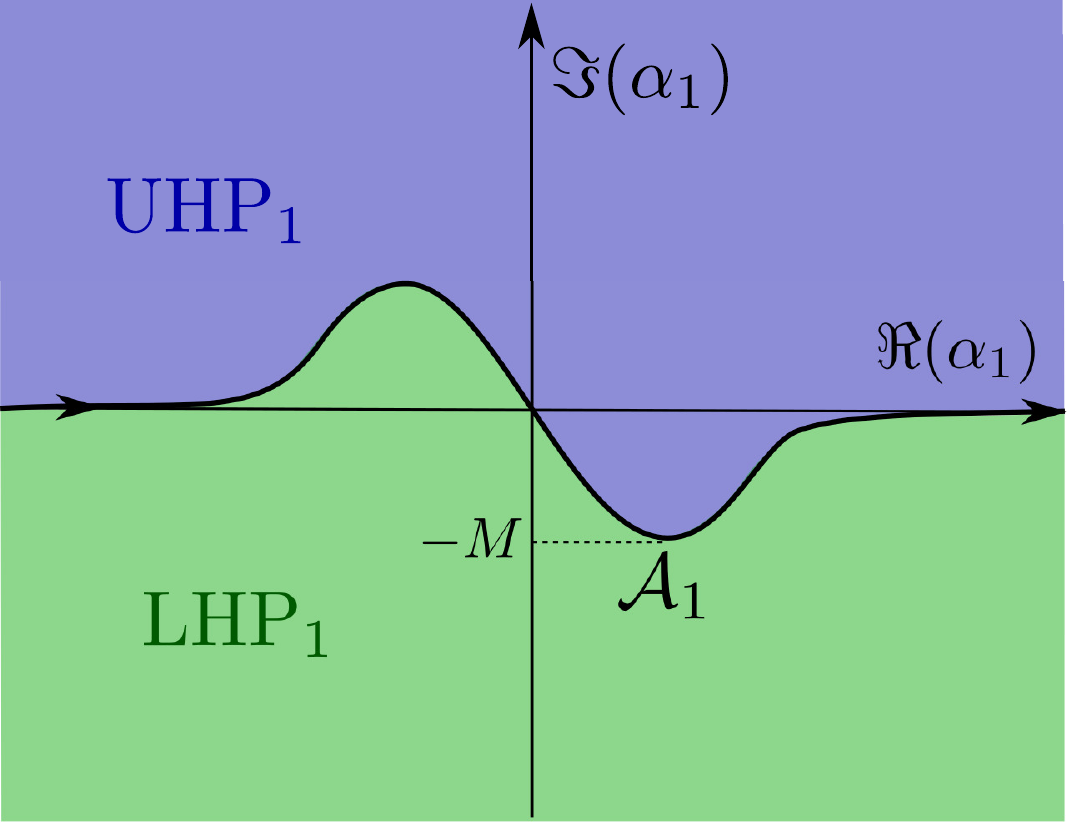}\quad\includegraphics[width=0.48\textwidth]{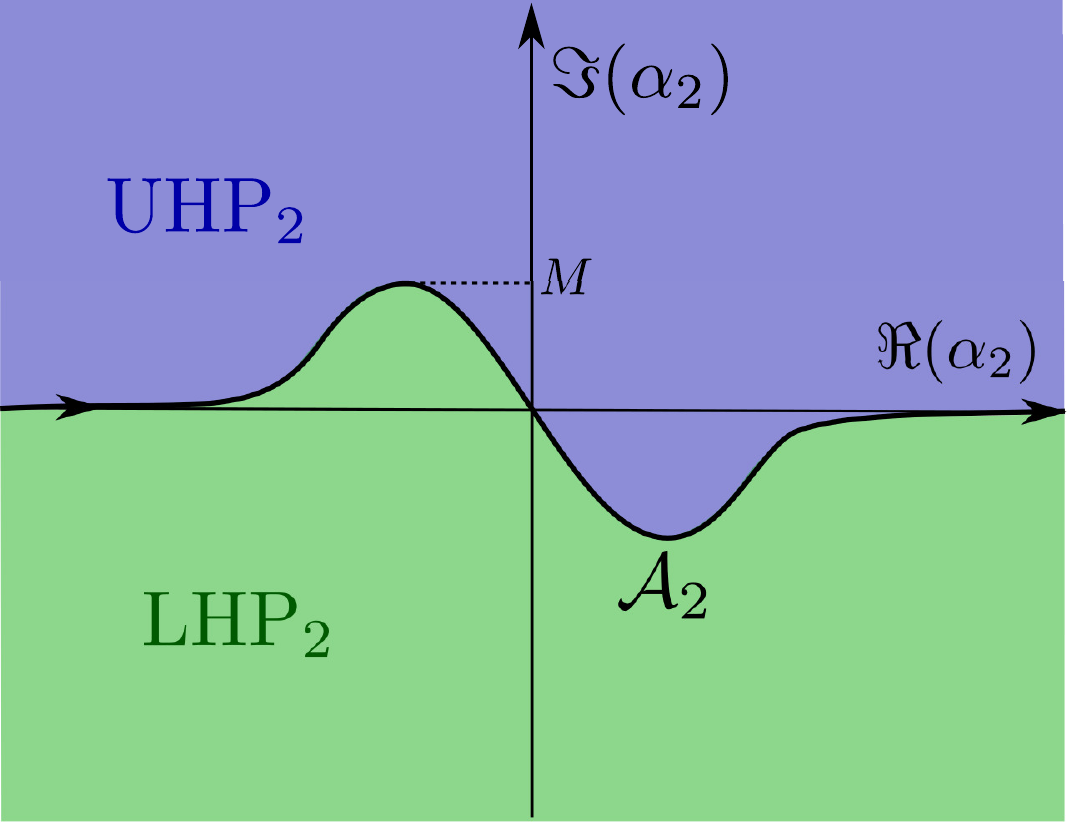}
  \caption{Description of the \RED{upper} and \RED{lower}-half planes $\tmop{UHP}_{1, 2}$
  and $\tmop{LHP}_{1, 2}$ lying on and above and on and below the integration contours
  $\mathcal{A}_{1, 2}$.}
\label{fig:contoursandhalfplanes}
\end{figure}

The quarter-plane is occupying the $(x_1 > 0, x_2 > 0, x_3 = 0)$ subspace of a
$(x_1, x_2, x_3)$ Cartesian space. The total physical wave field is denoted $u
(x_1, x_2, x_3)$ and the incident plane wave takes the form $e^{- i (a_1 x_1 +
a_2 x_2 + a_3 x_3)}$, where the constants $a_{1, 2, 3}$ depend solely on the
incident direction and the wave number $k > 0$ such that $a_1^2+a_2^2+a_3^2=k^2$. Note also that the time factor $e^{-i\omega t}$, where $\omega$ is the angular frequency of the wave, has been suppressed for brevity. This physical solution can be
expressed in the form of an inverse Fourier transform
\begin{eqnarray}
  u (\tmmathbf{x}, x_3) & = & \frac{1}{(2 \pi)^2} \int_{\mathcal{A}_1}
  \int_{\mathcal{A}_2} F_{+ +} (\tmmathbf{\alpha}) e^{- i\tmmathbf{\alpha}
  \cdot \tmmathbf{x}} e^{i \frac{x_3}{K (\tmmathbf{\alpha})}} \mathd \alpha_1
  \mathd \alpha_2,  \label{eq:inverseFouriersol}
\end{eqnarray}
where $\tmmathbf{x}= (x_1, x_2) \in \mathbb{R}^2$, $\tmmathbf{\alpha}=
(\alpha_1, \alpha_2) \in \mathbb{C}^2$ and the contours $\mathcal{A}_{1, 2}$
defined in Figure~\ref{fig:contoursandhalfplanes} naturally lead to the
definition of the upper and lower half planes $\tmop{UHP}_{1, 2}$ and
$\tmop{LHP}_{1, 2}$ of the $\alpha_{1, 2}$ complex planes. It is also useful
to define the domains $\mathcal{D}_{+ +} = \tmop{UHP}_1 \times \tmop{UHP}_2$
and $\mathcal{D}_{+ -} = \tmop{UHP}_1 \times \tmop{LHP}_2$. The crucial
unknown function $F_{+ +} (\tmmathbf{\alpha})$ is to be determined and $K
(\tmmathbf{\alpha}) = (k^2 - \alpha_1^2 - \alpha_2^2)^{- 1 / 2}$ is the
{\tmem{kernel}} of the problem.

One of the achievements of {\cite{Assier2019}} is the reduction of the
complicated problem of diffraction by a quarter-plane to two equations in the
two-complex-variables Fourier space. They relate the two unknown functions
$F_{+ +} (\tmmathbf{\alpha})$ and $G_{+ -} (\tmmathbf{\alpha})$, analytic on
$\mathcal{D}_{+ +}$ and $\mathcal{D}_{+ -}$ respectively, as follows
\begin{align}
  F_{+ +} &= F_{+ +}^{\star} +\mathcal{I} [G_{+ -}] \text{ \ on }
  \mathcal{D}_{+ +},  \label{eq:radlowoccurence}\\
  0 &= G_{+ -}^{\star} +\mathcal{J} [G_{+ -}]  \text{ on } \mathcal{D}_{+ -}, 
  \label{eq:compatibilityeq}
\end{align}
where $\mathcal{I}$ and $\mathcal{J}$ are explicitly-known Cauchy integral
operators depending on the kernel $K$, and the two functions $F^{\star}_{+ +}
(\tmmathbf{\alpha})$ and $G^{\star}_{+ -} (\tmmathbf{\alpha})$ are also known
explicitly and can be written in terms of $K$. \RED{It is remarkable that, in (\ref{eq:radlowoccurence}), the term $F_{++}^\star$ turns out to be exactly Radlow's erroneous ansatz \cite{Radlow1965}, while (\ref{eq:compatibilityeq}), the so-called \textit{compatibility equation}, only contains the unknown function $G_{+-}$. If (\ref{eq:compatibilityeq}) could somehow be inverted, then $G_{+-}$ would be known, and hence $F_{++}$ would be known by (\ref{eq:radlowoccurence}); see \cite{Assier2019} for a more detailed interpretation of these equations.}

In this article, we will focus on the function $G_{+ -}$, since it leads to
$F_{+ +}$ directly via (\ref{eq:radlowoccurence}) and hence to the sought-after $u$ via
(\ref{eq:inverseFouriersol}). Though unknown, this function can be expressed
in terms of the wave field $u$ by
\begin{eqnarray}
  G_{+ -} (\tmmathbf{\alpha}) & = & \int_{- \infty}^0 \int_0^{\infty} u
  (\tmmathbf{x}, 0) e^{i\tmmathbf{\alpha} \cdot \tmmathbf{x}} \mathd x_1
  \mathd x_2 .  \label{eq:integraldefofG+-}
\end{eqnarray}
The physical field $u$\RED{, restricted to the quadrant $(x_1\geq 0,x_2\leq 0,x_3=0)$,} must obey the following edge and vertex conditions:
\begin{eqnarray}
  u (\tmmathbf{x}, 0) & \underset{(x_1, x_2) \rightarrow (0^+, 0^-)}{\overset{|
  x_2 | \propto | x_1 |}{\approx}} & \mathcal{O} (r^{\RED{\nu} - 1 / 2}), \label{eq:1point6} \\
  u (\tmmathbf{x}, 0) & \underset{x_1 \rightarrow 0^+}{\overset{x_2 < 0 \text{
  } \tmop{fixed}}{\approx}} & \mathcal{O} (1), \label{eq:1point7} \\
  u (\tmmathbf{x}, 0) & \underset{x_2 \rightarrow 0^-}{\overset{x_1 > 0 \text{
  } \tmop{fixed}}{\approx}} & \mathcal{O} (| x_2 |^{1 / 2}), \label{eq:1point8}
\end{eqnarray}
where $r = \sqrt{x_1^2 + x_2^2}$ is the distance to the vertex and
\begin{eqnarray*}
  & \RED{\nu} = \sqrt{\RED{\ell}_1 + 1 / 4} \approx 0.7967, & 
\end{eqnarray*}
with $\RED{\ell}_1$ being the first eigenvalue of the Laplace-Beltrami operator (LBO)
with Dirichlet conditions on the cut defined by $\left\{ \theta =
\frac{\pi}{2}, \varphi \in \left[ 0, \frac{\pi}{2} \right] \right\}$ in the
usual spherical coordinates (see e.g. {\cite{Assier2012}},
{\cite{Assier2016}}). \RED{Physically, the restrictions (\ref{eq:1point6}) and (\ref{eq:1point8}) come from imposing that the energy remains bounded in the neighbourhood of the vertex $(x_1=0,x_2=0, x_3=0)$ and the edge $(x_1>0,x_2=0,x_3=0)$ respectively. The exact value of the exponent comes from a separation of variable argument in spherical coordinates for (\ref{eq:1point6}) and cylindrical polar coordinates (around the edge) for (\ref{eq:1point8}). These conditions are often referred to as Meixner conditions. The condition (\ref{eq:1point7}) is representing the fact that the field is perfectly well behaved since no element of the scatterer is encountered in this limit.}  

In the Fourier space, \RED{using the Abelian theorems (see e.g. \cite{Noble1958} ),} these conditions translate into the following asymptotic
behaviour for $G_{+ -}$:
\begin{eqnarray}
  G_{+ -} (\alpha_1, \alpha_2) & \underset{| \alpha_1 | \rightarrow
  \infty}{\overset{| \alpha_2 | \propto | \alpha_1 |}{=}} & \mathcal{O} \left(
  \frac{1}{| \alpha_1 |^{\RED{\nu} + 3 / 2}} \right),  \label{eq:asympI1}\\
  G_{+ -} (\alpha_1, \alpha_2) & \underset{| \alpha_1 | \rightarrow
  \infty}{\overset{\alpha_2 \tmop{fixed}}{=}} & \mathcal{O} \left( \frac{1}{|
  \alpha_1 |} \right),  \label{eq:asympI2}\\
  G_{+ -} (\alpha_1, \alpha_2) & \underset{| \alpha_2 | \rightarrow
  \infty}{\overset{\alpha_1 \tmop{fixed}}{=}} & \mathcal{O} \left( \frac{1}{|
  \alpha_2 |^{3 / 2}} \right),  \label{eq:asympI3}
\end{eqnarray}
for $\alpha_1 \in \tmop{UHP}_1$ and $\alpha_2 \in \tmop{LHP}_2$.


\RED{In order to derive the equations
(\ref{eq:radlowoccurence})--(\ref{eq:compatibilityeq}) rigorously in {\cite{Assier2019}}, we had to prove that $G_{+-}$ satisfies an important property, i.e.\ that $\mathcal{I} [G_{+ -}] (\alpha_1, \alpha_2)$
tends to zero as $| \alpha_2 | \rightarrow \infty$.}

\RED{The purpose of the present article is two-fold; first
to suggest an efficient approximation scheme to solve the quarter-plane
problem and second to highlight a new canonical special function. The former} will be done by introducing the explicitly defined integral $I_{+ -}^{\varepsilon}$, and by showing that it mimics the behaviour of $G_{+ -}$. By
this we mean that $I_{+ -}^{\varepsilon}$ should be analytic on
$\mathcal{D}_{+ -}$, should have the asymptotic behaviour
(\ref{eq:asympI1})--(\ref{eq:asympI3}) and should satisfy $\mathcal{I} [I_{+
-}^{\varepsilon}] \rightarrow 0$ as $| \alpha_2 | \rightarrow \infty$.

The rest of the paper is organised as follows: in Section
\ref{sec:studyofintegral}, we give the integral expression of $I_{+
-}^{\varepsilon}$ again and explain where it comes from; in Section
\ref{sec:eigenfunction} we highlight some important properties of the first
Laplace-Beltrami eigenfunction; in Section \ref{sec:asymptoticbehaviour} we
prove that $I_{+ -}^{\varepsilon}$ does indeed have the correct asymptotic
behaviour (\ref{eq:asympI1})--(\ref{eq:asympI3}); and in Section
\ref{sec:sec5}, we prove that $\mathcal{I} [I_{+ -}^{\varepsilon}]$ does tend
to zero as $| \alpha_2 |$ tends to infinity. Finally, we discuss the
implications of our findings and conclude the paper in Section
\ref{sec:conclusion}.

\section{\RED{A canonical} integral}\label{sec:studyofintegral}

In this section, we aim to derive and construct explicitly a function defined
by an integral that satisfies all the conditions required of $G_{+ -}$.
Starting from the integral representation (\ref{eq:integraldefofG+-}) and
using the change of variable $x_1 = r \cos (\varphi)$, $x_2 = r \sin
(\varphi)$, for $r \in \mathbb{R}^+$ and $\varphi \in [3 \pi / 2, 2 \pi]$, we
can write
\begin{eqnarray}
  G_{+ -} (\tmmathbf{\alpha}) & = & \int_{\frac{3 \pi}{2}}^{2 \pi}
  \int_0^{\infty} \hat{u} (r, \varphi) e^{ir (\alpha_1 \cos (\varphi) +
  \alpha_2 \sin (\varphi))} r \mathd r \mathd \varphi, 
  \label{eq:changevariableintegralG}
\end{eqnarray}
where $\hat{u} (r, \varphi) = u (r \cos (\varphi), r \sin (\varphi), 0)$.
Moreover, as explained in {\cite{Assier2019}} for example, using separation of
variables, it can be shown that
\begin{eqnarray}
  \hat{u} (r, \varphi) & \underset{r \rightarrow 0}{\overset{\varphi \in
  \left[ \frac{3 \pi}{2}, 2 \pi \right)}{\sim}} & A \RED{f} \left( \frac{\pi}{2},
  \varphi \right) r^{\RED{\nu} - 1 / 2},  \label{eq:uhatbehaviouratzero}
\end{eqnarray}
where $\RED{f} (\theta, \varphi)$ is the eigenfunction of the LBO associated to
the first eigenvalue $\lambda_1$ and $A$ is a constant. For technical reasons
that will become apparent later on, let us rewrite the asymptotic behaviour
(\ref{eq:uhatbehaviouratzero}) in a slightly different form
\begin{eqnarray}
  \hat{u} (r, \varphi) & \underset{r \rightarrow 0}{\sim} & A \RED{f} \left(
  \frac{\pi}{2}, \varphi \right) r^{\RED{\nu} - 1 / 2} e^{- \varepsilon r}, 
  \label{eq:newuhatbehaviour}
\end{eqnarray}
for some $\varepsilon > 0$. Note that (\ref{eq:uhatbehaviouratzero}) and
(\ref{eq:newuhatbehaviour}) are equivalent to leading order since $e^{-
\varepsilon r} \rightarrow 1$ as $r \rightarrow 0$. Because
\RED{the asymptotic behaviour as $|\alpha_{1,2}|\to\infty$}
in the Fourier space is intrinsically linked to the near-field
behaviour in the physical space, we are interested in the integral $I$
obtained by replacing $\hat{u} (r, \varphi)$ by its leading order behaviour
(\ref{eq:newuhatbehaviour}) in (\ref{eq:changevariableintegralG}):
\begin{eqnarray*}
  \PUR{I} & = &
A \int_{\frac{3 \pi}{2}}^{2 \pi} \RED{f} \left( \frac{\pi}{2}, \varphi
  \right) \int_0^{\infty} r^{\RED{\nu} + 1 / 2} e^{ir (\alpha_1 \cos (\varphi) +
  \alpha_2 \sin (\varphi))} e^{- \varepsilon r} \mathd r \mathd \varphi .
\end{eqnarray*}
Note that the integral over $r$ takes the form
\begin{align}
  \int_0^{\infty} r^{\RED{\lambda} - 1} e^{- \mu r} \mathd r,  
\end{align}
for $\RED{\lambda} = \RED{\nu} + 3 / 2$ and $\mu = - i (\alpha_1 \cos (\varphi) + \alpha_2
\sin (\varphi) + i \varepsilon)$. Note (see {\cite{GRADSHTEYN1980211}} p317,
3.381.4.) that for $\tmop{Re} (\RED{\lambda}) > 0$ and $\tmop{Re} (\mu) > 0$, this
integral is exactly equal to $\frac{\Gamma (\RED{\lambda})}{\mu^{\RED{\lambda}}}$, where $\Gamma$
is the Euler Gamma function. It is clear that $\tmop{Re} (\RED{\lambda}) > 0$ since
$\RED{\nu} \geqslant 0$. Moreover, we can choose $\varepsilon$ such that
$\tmop{Re} (\mu) > 0$. In order to do so, we refer to Figure
\ref{fig:contoursandhalfplanes}, to see that for all $\tmmathbf{\alpha} \in
\mathcal{D}_{+ -}$, we have $\tmop{Im} (\alpha_1) \geqslant - M$ and
$\tmop{Im} (\alpha_2) \leqslant M$ for some $M > 0$ depending on the choice of
contour $\mathcal{A}_{1, 2}$. Remembering that for $\varphi \in \left[ \frac{3
\pi}{2}, 2 \pi \right]$ we have $0 \leqslant \cos (\varphi) \leqslant 1$ and
$- 1 \leqslant \sin (\varphi) \leqslant 0$, it is possible to show that upon
choosing $\varepsilon$ such that $\varepsilon > 2 M$, we have $\tmop{Re} (\mu)
> 0$ for all $\tmmathbf{\alpha} \in \mathcal{D}_{+ -}$, and hence $I$ can be
rewritten as
\begin{eqnarray*}
  I & = & A \Gamma (\RED{\nu} + 3 / 2) \int_{\frac{3 \pi}{2}}^{2 \pi} \frac{\RED{f}
  \left( \frac{\pi}{2}, \varphi \right)}{(- i (\alpha_1 \cos (\varphi) +
  \alpha_2 \sin (\varphi) + i \varepsilon))^{\RED{\nu} + 3 / 2}} \mathd \varphi .
\end{eqnarray*}
Because of our choice of $\varepsilon$ the denominator of the integrand is
never zero for $\tmmathbf{\alpha} \in \mathcal{D}_{+ -}$; hence the
integral is a `$+ -$' function, i.e.\ it is analytic in $\mathcal{D}_{+ -}$. This naturally leads to the definition of the
\RED{canonical} integral $I_{+ -}^{\varepsilon}$ to be studied in this paper:
\begin{eqnarray}
  I_{+ -}^{\varepsilon} (\alpha_1, \alpha_2) & = & \int_{\frac{3 \pi}{2}}^{2
  \pi} \frac{\RED{f} \left( \frac{\pi}{2}, \varphi \right)}{(- i (\alpha_1 \cos
  (\varphi) + \alpha_2 \sin (\varphi) + i \varepsilon))^{\RED{\nu} + 3 / 2}}
  \mathd \varphi .  \label{eq:interestingintegral}
\end{eqnarray}
For the purpose of the present work, the values of $k$ and $M$ (and hence
$\varepsilon$) can be any strictly positive numbers. For numerical
illustration of our theoretical results, we will choose $k = 3$ and
$\varepsilon = 1$.

\section{A note on the Laplace-Beltrami
eigenfunction}\label{sec:eigenfunction}

Before deriving various properties of the newly-introduced integral $I_{+
-}^{\varepsilon}$, it is important to know the behaviour of the eigenfunction
$\RED{f} \left( \frac{\pi}{2}, \varphi \right)$. Its important properties are
summarised in the following lemma, the proof of which (linked to the physical
edge conditions) is omitted here for brevity.

\begin{lemma}
  \label{lem:eigenfunctionproperties}Let $\RED{f} (\theta, \varphi)$ be the first
  eigenfunction of the LBO. Then $\RED{f} \left( \frac{\pi}{2}, \varphi \right)$
  is equal to zero on the cut $\varphi \in \left[ 0, \frac{\pi}{2} \right]$,
  and is smooth and such that $0 < \RED{f} \left( \frac{\pi}{2}, \varphi \right)
  \leqslant 1$ for $\varphi \in \left( \frac{\pi}{2}, 2 \pi \right)$.
  Moreover, its behaviour at the edge of the non-zero region is given by
  \begin{eqnarray*}
    \RED{f} \left( \frac{\pi}{2}, \varphi \right) \underset{\varphi \rightarrow
    \pi / 2}{\overset{\varphi > \pi / 2}{=}} \mathcal{O} \left( \sqrt{\varphi
    - \pi / 2} \right) & \text{ and } & \RED{f} \left( \frac{\pi}{2}, \varphi
    \right) \underset{\varphi \rightarrow 2 \pi}{\overset{\varphi < 2 \pi}{=}}
    \mathcal{O} \left( \sqrt{2 \pi - \varphi} \right) .
  \end{eqnarray*}
  In particular, there exists a constant {\tmem{$\beta$}}, such that
  \begin{eqnarray*}
    \RED{f} \left( \frac{\pi}{2}, 2 \pi - \psi \right) & \underset{\psi
    \rightarrow 0^{}}{\overset{\psi > 0}{\sim}} & \beta \sqrt{\psi} .
  \end{eqnarray*}
  In addition, it transpires that $\RED{f} \left( \frac{\pi}{2}, \varphi \right)$
  is strictly decreasing for $\varphi \in \left[ \frac{3 \pi}{2}, 2 \pi
  \right]$.
\end{lemma}

\begin{figure}[h]
\centering
\includegraphics[width=0.6\textwidth]{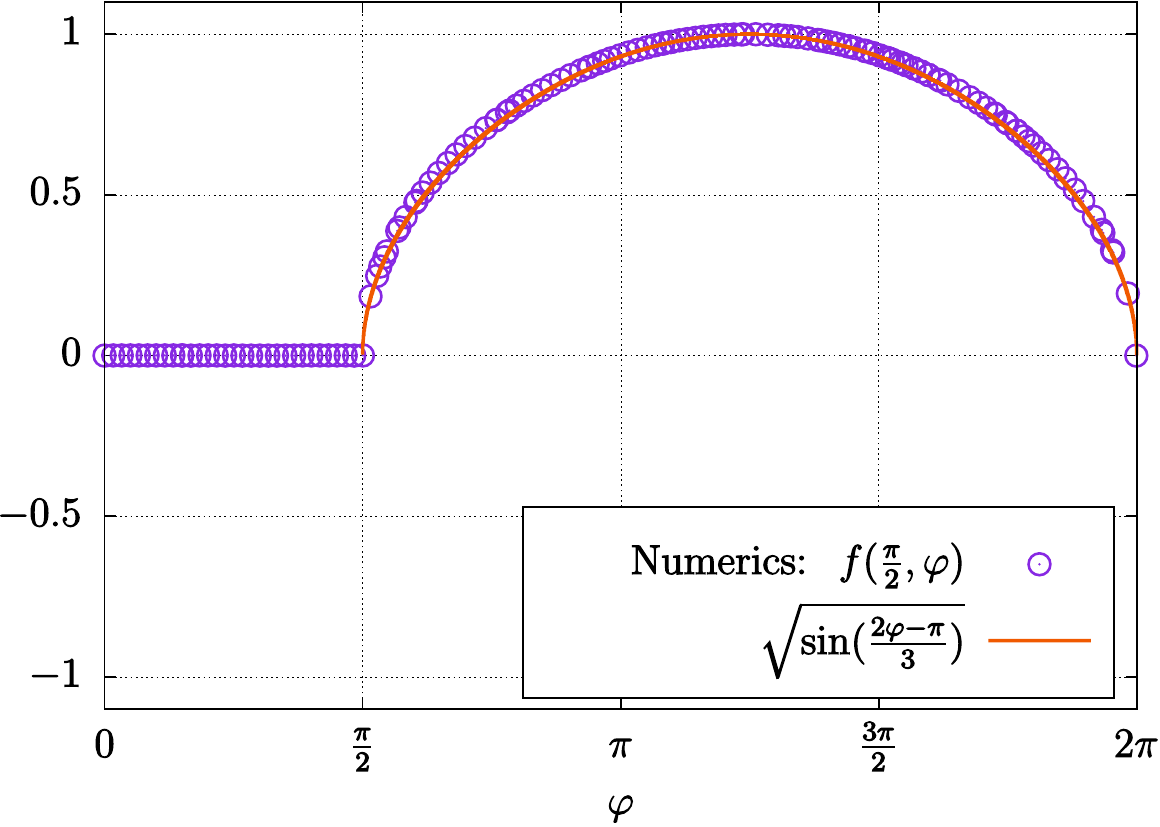}
  \caption{Comparison between a numerical approximation of $\RED{f} \left(
  \frac{\pi}{2}, \varphi \right)$ and the function $g (\varphi)$.}
\label{fig:linkeigenfunctionandsin}
\end{figure}

Though not an exact result, it seems that $\RED{f} \left( \frac{\pi}{2}, \varphi \right)$ is
approximated very well by the function $g (\varphi)$ defined by
\begin{eqnarray*}
  g (\varphi) & = & \left\{ \begin{array}{lll}
    0 & \tmop{if} & \varphi \in \left[ 0, \frac{\pi}{2} \right],\\
    \sqrt{\sin \left( \frac{2 \varphi - \pi}{3} \right)} & \tmop{if} & \varphi
    \in \left[ \frac{\pi}{2}, 2 \pi \right].
  \end{array} \right. 
\end{eqnarray*}
Note that $g$ trivially satisfies the conclusions of Lemma
\ref{lem:eigenfunctionproperties}, with $\beta = \sqrt{\frac{2}{3}}$. In
Figure~\ref{fig:linkeigenfunctionandsin}, we compare numerical results
obtained using a surface finite element method developed in
{\cite{Assier2016}} and the function $g$, showing excellent agreement.

\section{Asymptotic behaviour of $I_{+
-}^{\varepsilon}$}\label{sec:asymptoticbehaviour}

In this section, we show that the integral $I_{+ -}^{\varepsilon}
(\alpha_1, \alpha_2)$ shares the same rich asymptotic behaviour as $G_{+ -}
(\alpha_1, \alpha_2)$; that is, it should behave like
(\ref{eq:asympI1})--(\ref{eq:asympI3}), as $| \alpha_{1, 2} | \rightarrow
\infty$ within $\mathcal{D}_{+ -}$. We will summarise the key
results here, the detailed proofs being given in Appendix~\ref{app:app1}. Since
$\tmmathbf{\alpha}= (\alpha_1, \alpha_2) \in \mathcal{D}_{+ -}$, whenever $| \alpha_1 |$ (resp. $| \alpha_2 |$) tends
to infinity, we can write $\alpha_1 = | \alpha_1 | e^{i \phi_1}$ (resp.
$\alpha_2 = | \alpha_2 | e^{i \phi_2}$) for $\phi_1 \in (0, \pi)$ (resp.
$\phi_2 \in (- \pi, 0)$). Note that if $\alpha_{1,2}$ are not assumed large, we cannot write them in this way, since they may lie within the indented part of the contours $\mathcal{A}_{1,2}$.

\paragraph{Asymptotic behaviour when both $| \alpha_1 |$ and $| \alpha_2 |$
tend to infinity within $\mathcal{D}_{+ -}$}

\begin{figure}[h]
\centering
 \includegraphics[width=0.49\textwidth]{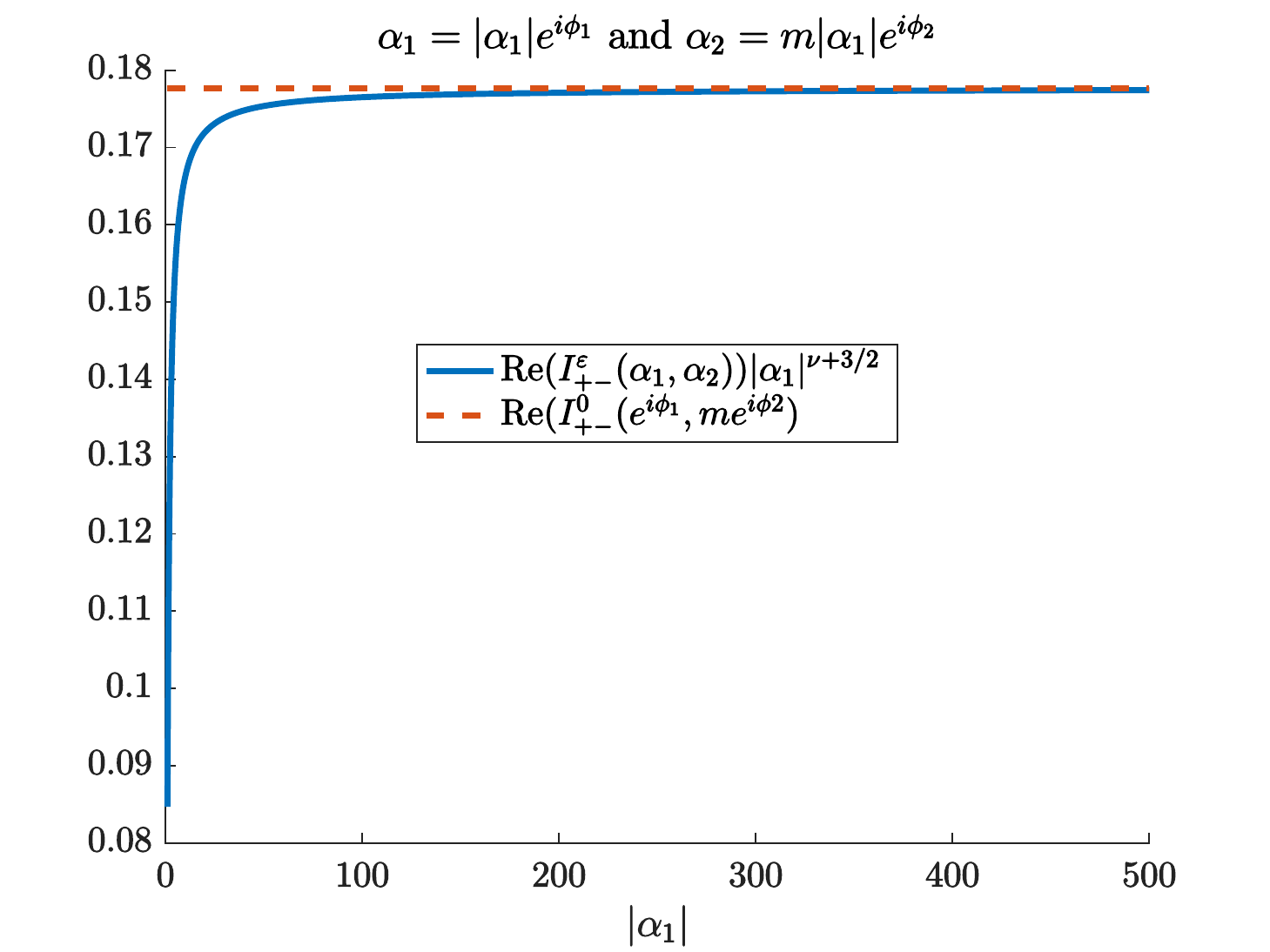}
\includegraphics[width=0.49\textwidth]{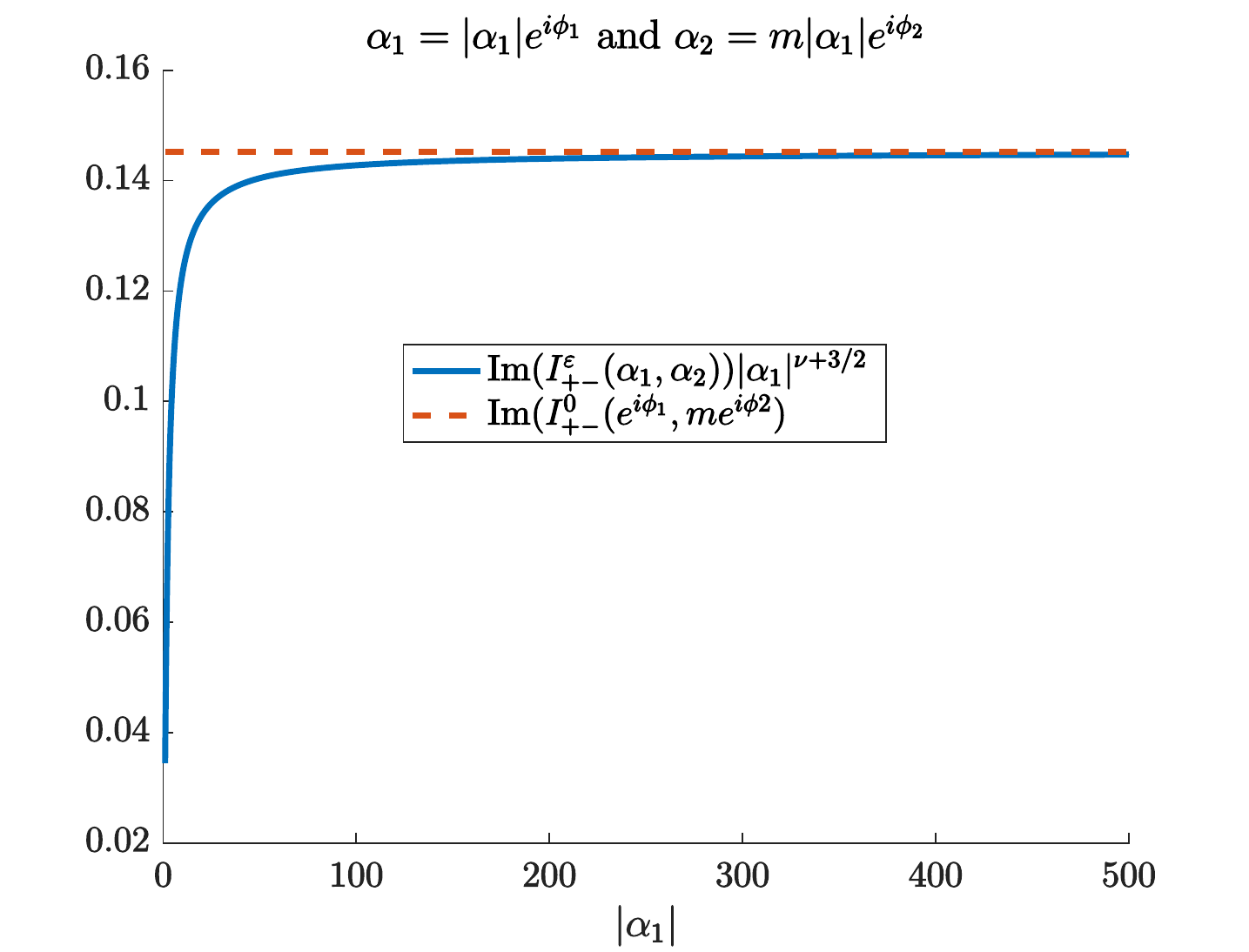}
  \caption{Numerical illustration of the asymptotic behaviour
  (\ref{eq:asymp-behav-1}) as both $| \alpha_{1, 2} | \rightarrow \infty$ for
  $\phi_1 = \frac{\pi}{4}$, $\phi_2 = - \frac{\pi}{2}$, $\varepsilon = 1$ and
  $m = 2$, using $g (\varphi)$ instead of $\RED{f} \left( \frac{\pi}{2}, \varphi
  \right)$ in the definition of $I_{+ -}^{\varepsilon}$.}
\label{fig:asymp1}
\end{figure}

This is the simplest of the three different cases to be considered. We will
take both $| \alpha_1 |$ and $| \alpha_2 | \rightarrow \infty$ within
$\mathcal{D}_{+ -}$, in such a way that there exists an $m > 0$ such that $|
\alpha_2 | = m | \alpha_1 |$. In this case, we can write $\alpha_{1, 2} = |
\alpha_{1, 2} | e^{i \phi_{1, 2}}$ and we have
\begin{eqnarray}
  I_{+ -}^{\varepsilon} (\alpha_1, \alpha_2) & \underset{| \alpha_1 |
  \rightarrow \infty}{\overset{| \alpha_2 | = m | \alpha_1 |}{\sim}} &
  \frac{1}{| \alpha_1 |^{\RED{\nu} + 3 / 2}} I_{+ -}^0 (e^{i \phi_1}, me^{i
  \phi_2}) .  \label{eq:asymp-behav-1}
\end{eqnarray}
The validity of this asymptotic behaviour is illustrated in Figure
\ref{fig:asymp1}.

\paragraph{Asymptotic behaviour when $| \alpha_1 | \rightarrow \infty$ within
$\tmop{UHP}_1$ and $\alpha_2$ is fixed in $\tmop{LHP}_2$}

\begin{figure}[h]
\centering
  \includegraphics[width=0.49\textwidth]{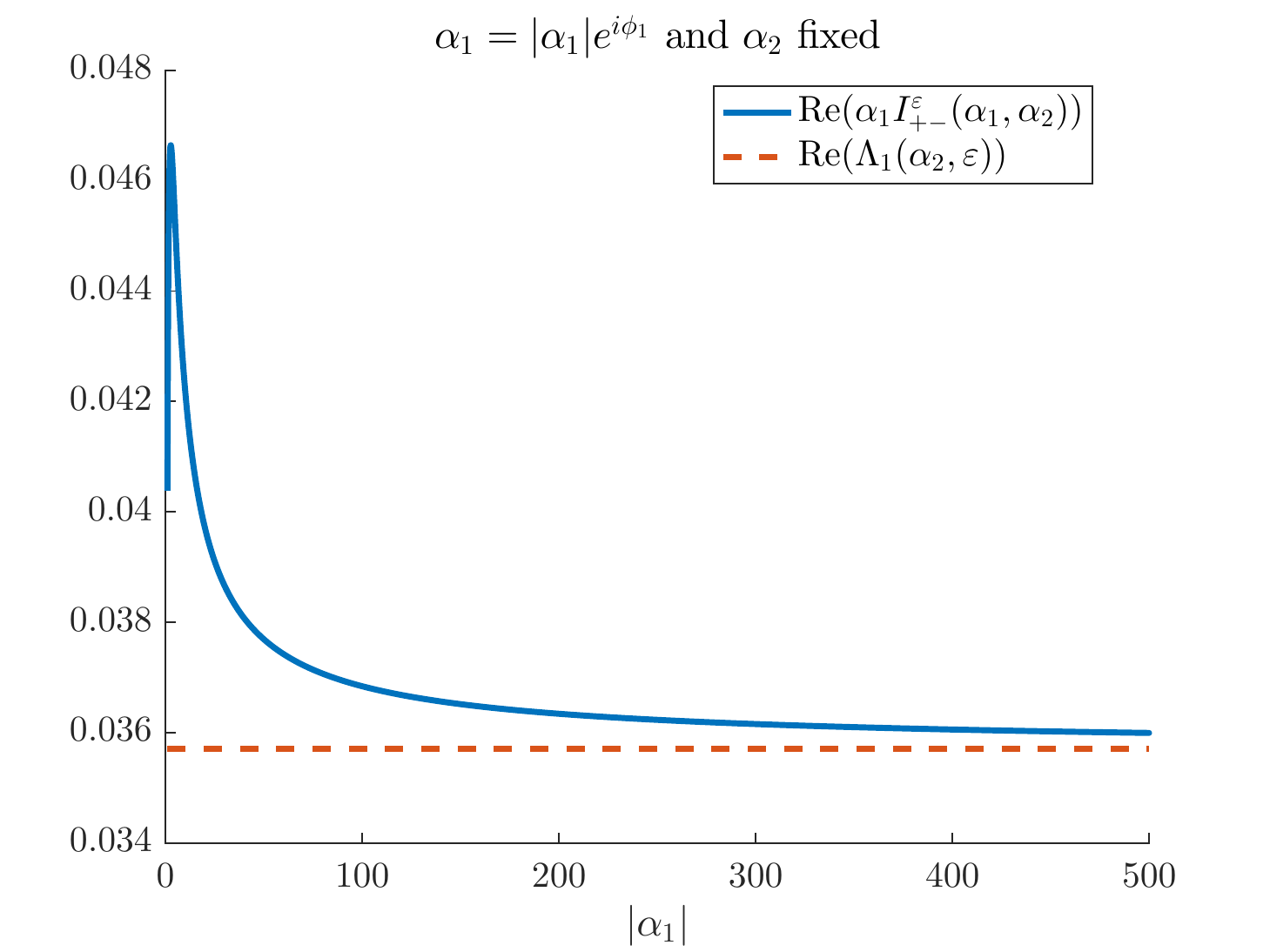}
\includegraphics[width=0.49\textwidth]{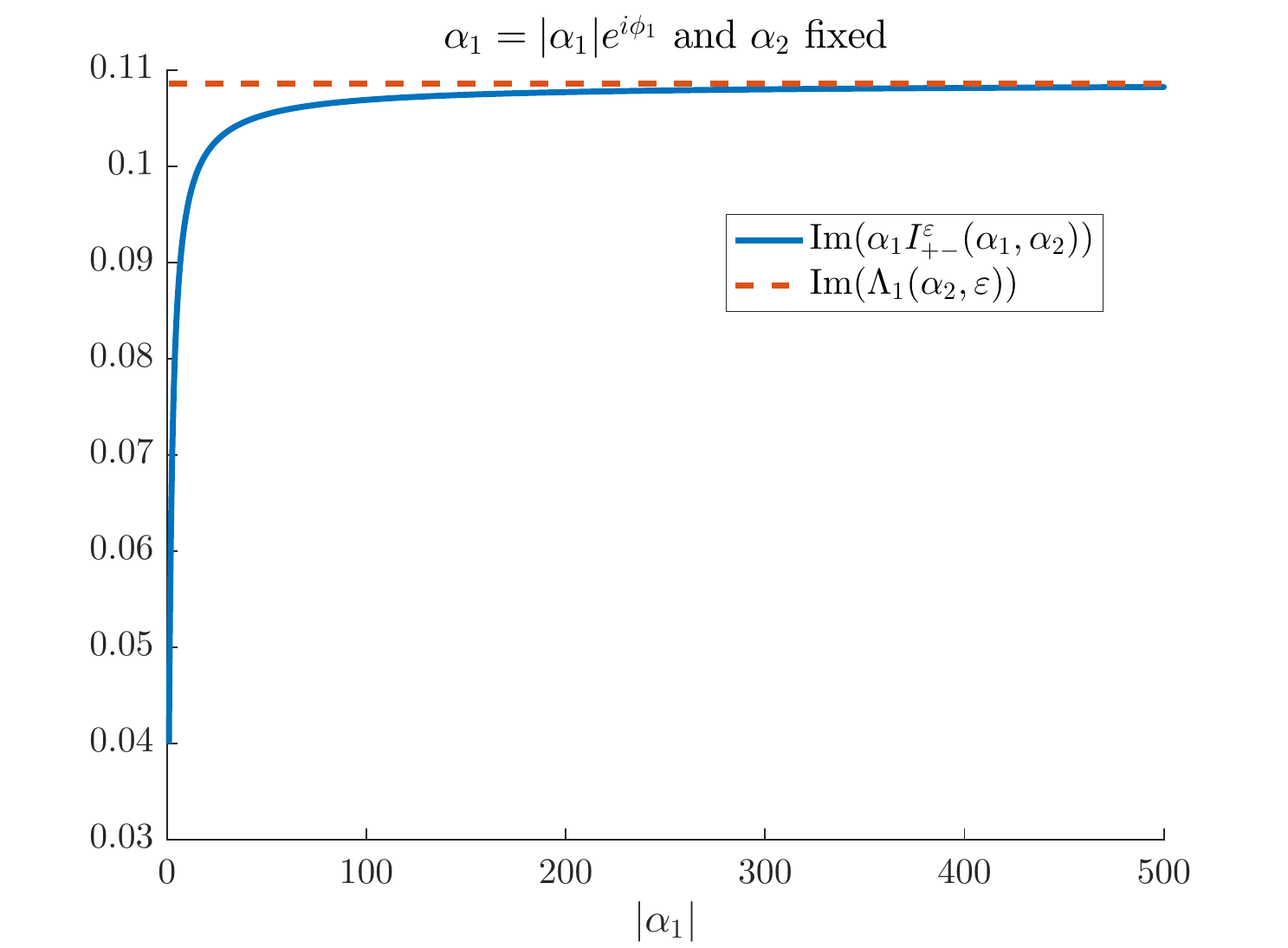}
  \caption{Numerical illustration of the asymptotic behaviour
  (\ref{eq:asymp-behav-2}) as $| \alpha_1 | \rightarrow \infty$ for $\phi_1 =
  \frac{\pi}{4}$, fixed $\alpha_2 = 1 - 3 i$ and $\varepsilon = 1$, using $g
  (\varphi)$ instead of $\RED{f} \left( \frac{\pi}{2}, \varphi \right)$ in the
  definition of $I_{+ -}^{\varepsilon}$.}
\label{fig:asymp2}
\end{figure}

In this case, we can write $\alpha_1 = | \alpha_1 | e^{i \phi_1}$, and we have
\begin{eqnarray}
  I_{+ -}^{\varepsilon} (\alpha_1, \alpha_2) & \underset{| \alpha_1 |
  \rightarrow \infty}{\overset{\alpha_2 \tmop{fixed}}{\sim}} & \frac{\Lambda_1
  (\alpha_2, \varepsilon)}{\alpha_1},  \label{eq:asymp-behav-2}
\end{eqnarray}
where
\begin{eqnarray*}
  \Lambda_1 (\alpha_2, \varepsilon) & = & \frac{2 i\RED{f} \left( \frac{\pi}{2},
  \frac{3 \pi}{2} \right)}{(1 + 2 \RED{\nu}) } \times \frac{1}{(\varepsilon + i
  \alpha_2)^{\RED{\nu} + 1 / 2}} \cdot
\end{eqnarray*}
As can be seen in Appendix~\ref{app:app1}, the proof is slightly more subtle
than the previous case, and one needs to split the $\varphi$ integral into two
parts, one where $\cos (\varphi)$ is very small, and one where it is bounded
away from zero. The validity of the asymptotic behaviour
(\ref{eq:asymp-behav-2}) is illustrated in Figure~\ref{fig:asymp2}.

\paragraph{Asymptotic behaviour when $| \alpha_2 | \rightarrow \infty$ within
$\tmop{LHP}_2$ and $\alpha_1$ is fixed in $\tmop{UHP}_1$}

\begin{figure}[h]
\centering
  \includegraphics[width=0.5\textwidth]{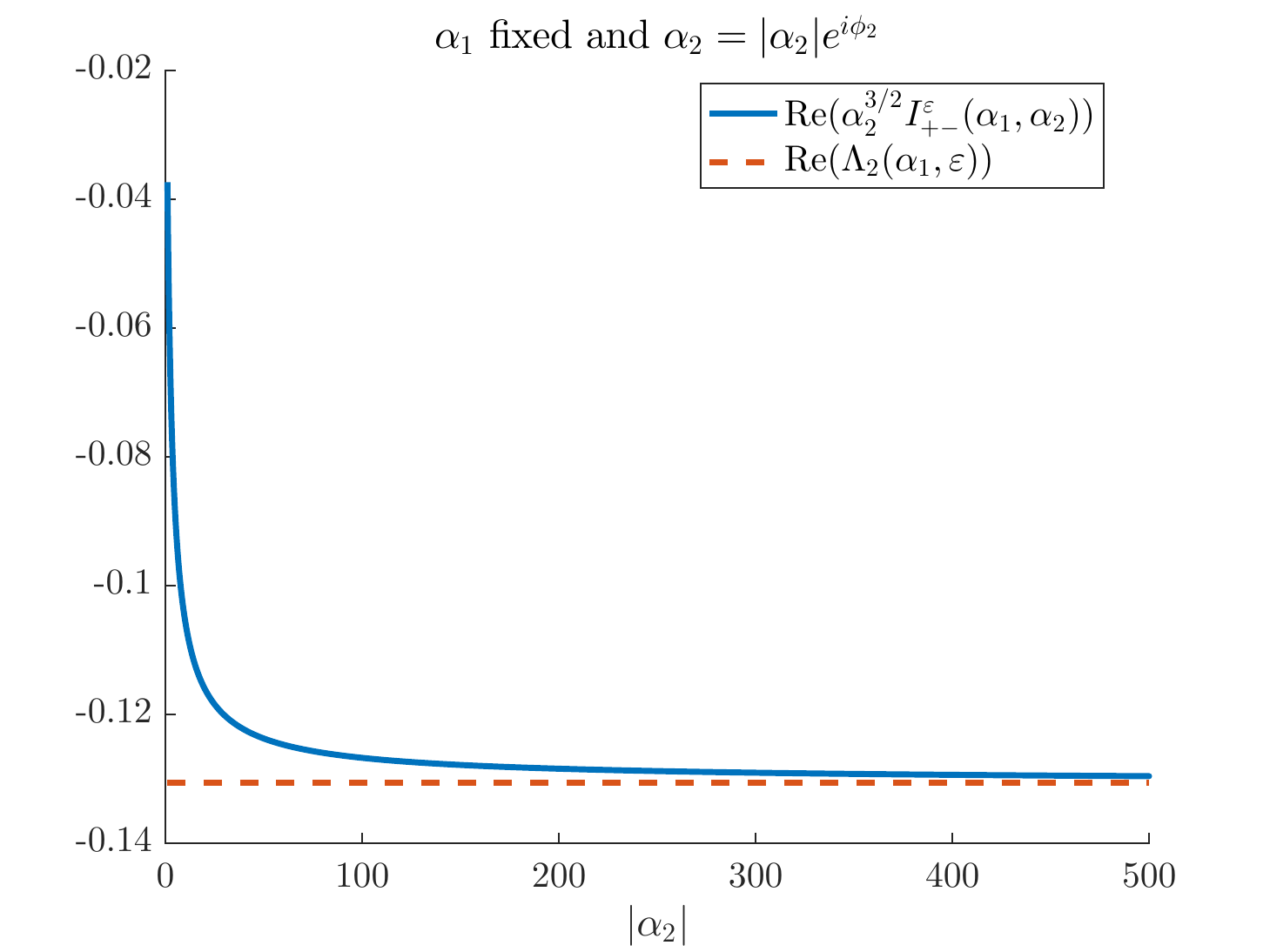}\includegraphics[width=0.5\textwidth]{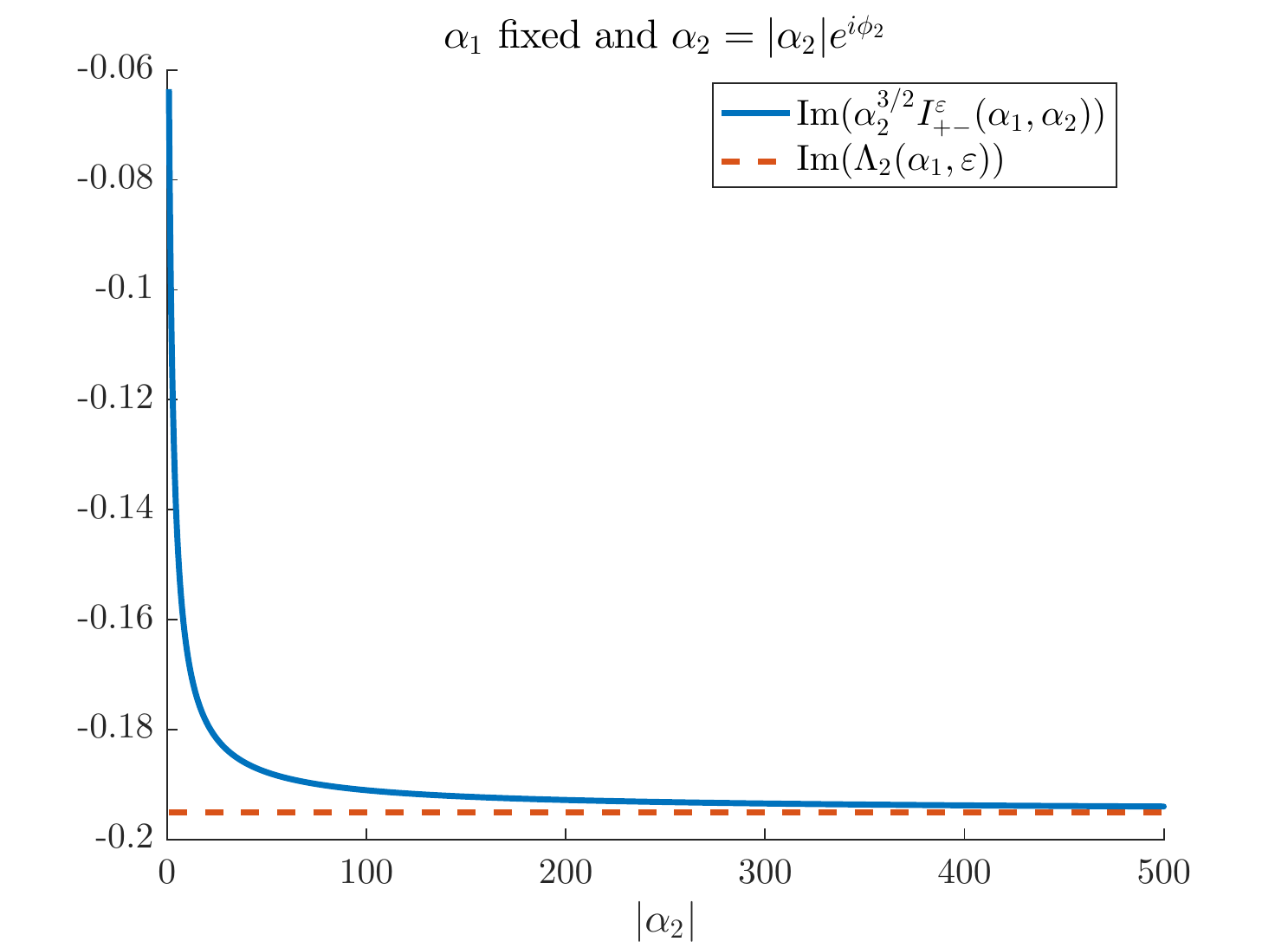}
  \caption{Numerical illustration of the asymptotic behaviour of
  (\ref{eq:asymp-behav-3}) as $| \alpha_2 | \rightarrow \infty$ for $\phi_2 =
  - \frac{\pi}{4}$, fixed $\alpha_1 = 1 + 3 i$ and $\varepsilon = 1$, using $g
  (\varphi)$ instead of $\RED{f} \left( \frac{\pi}{2}, \varphi \right)$ in the
  definition of $I_{+ -}^{\varepsilon}$.}
\label{fig:asymp3}
\end{figure}

In this case, we can write $\alpha_2 = | \alpha_2 | e^{i \phi_2}$, and we have
\begin{eqnarray}
  I_{+ -}^{\varepsilon} (\alpha_1, \alpha_2) & \underset{| \alpha_2 |
  \rightarrow \infty}{\overset{\alpha_1 \tmop{fixed}}{\sim}} & \frac{\Lambda_2
  (\alpha_1, \varepsilon)}{\alpha_2^{3 / 2}},  \label{eq:asymp-behav-3}
\end{eqnarray}
where
\begin{eqnarray*}
  \Lambda_2 (\alpha_1, \varepsilon) & = & \frac{\beta \sqrt{\pi} \Gamma
  (\RED{\nu}) e^{-3i\frac{\pi}{4}}}{2 \Gamma (\RED{\nu} + 3 / 2)} \times
  \frac{1}{(\varepsilon - i \alpha_1)^{\RED{\nu}}} \cdot
\end{eqnarray*}
Here again, as can be seen in Appendix~\ref{app:app1}, the proof is subtle and
requires a particular split of the $\varphi$ integral, this time according to
whether or not $\sin (\varphi)$ is close to being zero. The validity of the
asymptotic behaviour (\ref{eq:asymp-behav-3}) is illustrated in Figure
\ref{fig:asymp3}.

\begin{remark}
  Even though the results in this section have been derived for $\phi_1 \in
  (0, \pi)$ and $\phi_2 \in (- \pi, 0)$, they remain valid for $\phi_1 = 0,
  \pi$ and $\phi_2 = - \pi, 0$; this is easily checked numerically. Hence these
  asymptotic results can be used when $| \alpha_{1, 2} | \rightarrow \infty$
  along $\mathcal{A}_{1, 2}$.
\end{remark}

\begin{remark}
  If we let $| \alpha_2 | \rightarrow \infty$ in (\ref{eq:asymp-behav-2}), we
  obtain a quantity behaving like $\mathcal{O} (\alpha_1^{- 1} /
  \alpha_2^{\RED{\nu} + 1 / 2}),$while if we let $| \alpha_1 | \rightarrow \infty$
  in (\ref{eq:asymp-behav-3}), we obtain a quantity behaving like $\mathcal{O}
  (\alpha_1^{- \RED{\nu}} / \alpha_2^{3 / 2})$, both expressions being compatible
  with the behaviour (\ref{eq:asymp-behav-1}), for $|\alpha_1| \propto  |\alpha_2|$.
\end{remark}

We have hence shown that $I_{+ -}^{\varepsilon}$ and $G_{+ -}$ have the same
asymptotic behaviour at infinity. 
\RED{We mentioned in the introduction that in order to prove the main result of
{\cite{Assier2019}}, that is (\ref{eq:radlowoccurence})--(\ref{eq:compatibilityeq}), it was crucial to show that
$\mathcal{I} [G_{+ -}]\to 0$ as $| \alpha_2 | \rightarrow \infty$. We hence expect that a good approximation to $G_{+-}$ should have the same property.}
In the following
section we will precisely \RED{define the integral operator $\mathcal{I}$ and show that $I_{+ -}^{\varepsilon}$ does indeed satisfy this crucial property.}


\section{On the behaviour of $\mathcal{I} [I_{+ -}^{\varepsilon}]$ as $|
\alpha_2 | \rightarrow \infty$}\label{sec:sec5}

\subsection{The integral operator $\mathcal{I}$}

\RED{For \RED{a generic} function $\Phi (\alpha_1, \alpha_2)$ \RED{analytic on $\mathcal{A}_1\times\mathcal{A}_2$}, the Cauchy integral
operator $\mathcal{I}$ is defined as follows:
\begin{eqnarray}
  \mathcal{I} [\Phi] & = & \left[ \frac{K_{- +} (a_1, \alpha_2)}{K_{+ -}
  (\tmmathbf{\alpha})} \left[ \frac{\Phi}{K_{- \circ}} \right]_{+ \circ}
  \right]_{\circ +}, 
\end{eqnarray}
where \RED{the functions $K_{-+}$, $K_{+-}$ and} $K_{- \circ} (\tmmathbf{\alpha}) = K_{- -} (\tmmathbf{\alpha}) K_{- +}
(\tmmathbf{\alpha})$ \RED{are directly related to}}
the
four-way factorisation of the kernel $K (\alpha_1, \alpha_2)$ discussed in
{\cite{Assier2019}}. On $\mathcal{A}_1 \times \mathcal{A}_2$, the kernel $K$
can be written as
\begin{eqnarray*}
  K (\tmmathbf{\alpha}) & = & K_{+ +} (\tmmathbf{\alpha}) K_{+ -}
  (\tmmathbf{\alpha}) K_{- +} (\tmmathbf{\alpha}) K_{- -} (\tmmathbf{\alpha}),
\end{eqnarray*}
where $K_{++} (\tmmathbf{\alpha})$ is analytic on $\mathcal{D}_{++}$, etc. Explicit integral expressions for these factors, which may be evaluated very rapidly, are given in
{\cite{Assier2019}}, and the following asymptotic behaviour is valid
\begin{eqnarray}
  K_{\pm \pm} (\alpha_1, \alpha_2) & \overset{\alpha_1 \text{
  fixed}}{\underset{| \alpha_2 | \rightarrow \infty}{=}} & \mathcal{O} (1 / |
  \alpha_2 |^{1 / 4}) \text{ for } \tmmathbf{\alpha} \in \mathcal{D}_{\pm \pm}
  .  \label{eq:alpha2behaviourK++}
\end{eqnarray}
\RED{T}he brackets $\left[ \hspace{1em} \right]_{+ \circ}$
and $\left[ \hspace{1em} \right]_{\circ +}$ 
are Cauchy integral sum-split
operators in the $\alpha_1$ and $\alpha_2$ complex planes respectively,
defined 
for a generic function $\Phi$ 
by
\begin{eqnarray*}
  {}[\Phi]_{+ \circ} (\alpha_1, \alpha_2) = \frac{1}{2\pi i}
  \int_{\mathcal{A}^b} \frac{\Phi (z, \alpha_2)}{(z - \alpha_1)} \mathd z &
  \text{ and } & [\Phi]_{\circ +} (\alpha_1, \alpha_2) = \frac{1}{2 \pi i}
  \int_{\mathcal{A}^b} \frac{\Phi (\alpha_1, z)}{(z - \alpha_2)} \mathd z,
\end{eqnarray*}
where $\mathcal{A}^b$ is a contour that lies just below $\mathcal{A}_{1}$ or $\mathcal{A}_{2}$ as appropriate.
This ensures that $[\Phi]_{+ \circ}$ and $[\Phi]_{\circ +}$ can be freely
evaluated (and are analytic) on $\mathcal{A}_{1, 2}$.

As discussed for example in \cite{Assier2019} and
{\cite{Assier2018a}}, the function \RED{$K_{-\circ}$} can be written analytically from the
$\alpha_1$ factorisation of $K (\tmmathbf{\alpha})$, and is
\begin{eqnarray}
  K_{- \circ} (\alpha_1, \alpha_2) & = & \frac{1}{\sqrt{\sqrt{k^2 -
  \alpha_2^2} - \alpha_1}},  \label{eq:defKminuso}
\end{eqnarray}
with a careful choice of branch-cut location (in $\text{UHP}_1$ when the function is seen as a function of $\alpha_1$). We also remind the reader that $a_1$ is a constant depending on the incident angles.

\RED{Now that we have defined properly what was meant by $\mathcal{I}[I_{+
-}^{\varepsilon}]$, we will prove in the two next sections that it does indeed tend to zero as $|\alpha_2|$ tends to infinity.}

\subsection{A sufficient condition}

Because of (\ref{eq:alpha2behaviourK++}), it is clear that
\begin{eqnarray}
  \frac{K_{- +} (a_1, \alpha_2)}{K_{+ -} (\tmmathbf{\alpha})} & \underset{|
  \alpha_2 | \rightarrow \infty}{\overset{\alpha_1 \in \tmop{UHP}_1}{=}} &
  \mathcal{O} (1) .  \label{eq:cond-to-tend-to-one}
\end{eqnarray}
It is well-known (see e.g. Lemma B.1 of {\cite{Assier2019}}) that if $\Phi
(\alpha_1, \alpha_2)$ tends to zero as a power of $| \alpha_2 |$ as $|
\alpha_2 | \rightarrow \infty$ on $\mathcal{A}_2$, then the sum-split bracket
$[\Phi]_{\circ +}$ tends to zero as $| \alpha_2 | \rightarrow \infty$. Hence
for $\mathcal{I} [I_{+ -}^{\varepsilon}]$ to tend to zero as $| \alpha_2 |
\rightarrow \infty$, using (\ref{eq:cond-to-tend-to-one}), it is enough to
show that
\begin{eqnarray}
  \left[ \frac{I_{+ -}^{\varepsilon}}{K_{- \circ}} \right]_{+ \circ}
  \!\!\!\!(\alpha_1, \alpha_2) & \underset{| \alpha_2 | \rightarrow
  \infty}{\overset{\alpha_2 \in \mathcal{A}_2}{=}} & \mathcal{O} (1 / |
  \alpha_2 |^{\gamma}),  \label{eq:condition}
\end{eqnarray}
for some $\gamma > 0$, which remains a non-trivial task that will be completed in what
follows.

\subsection{Proof strategy }\label{sec:proofstrategy}

\begin{figure}[h]
\centering
  \includegraphics[width=0.49\textwidth]{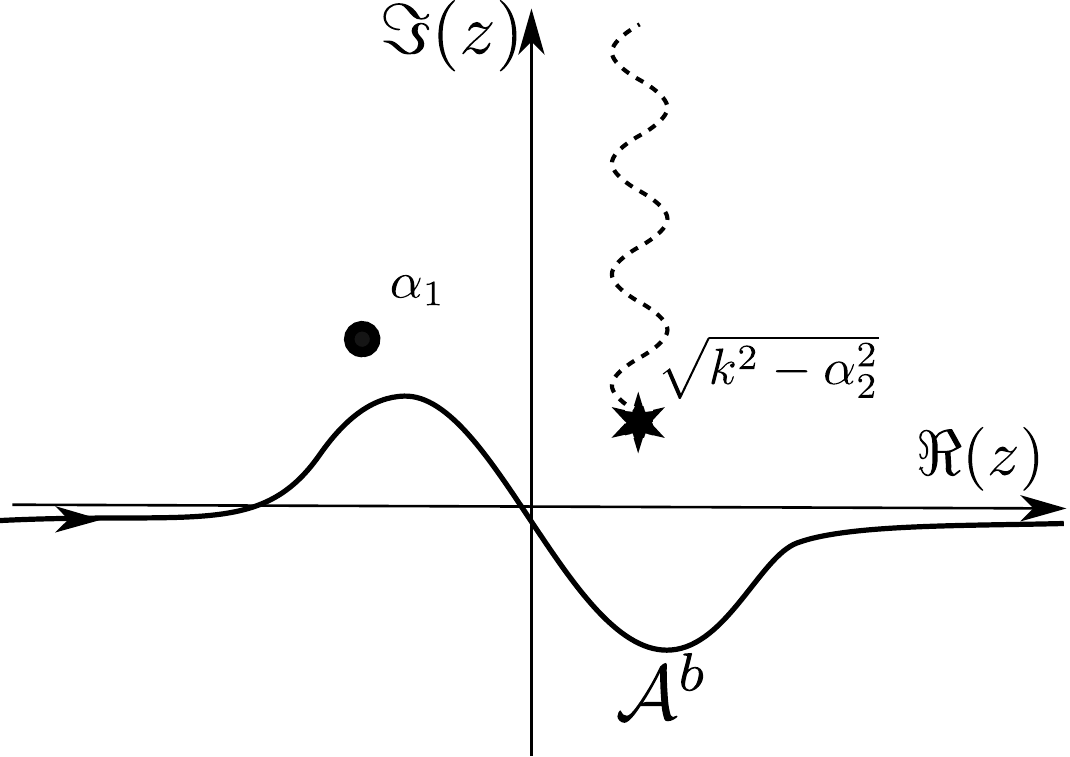}
 \includegraphics[width=0.49\textwidth]{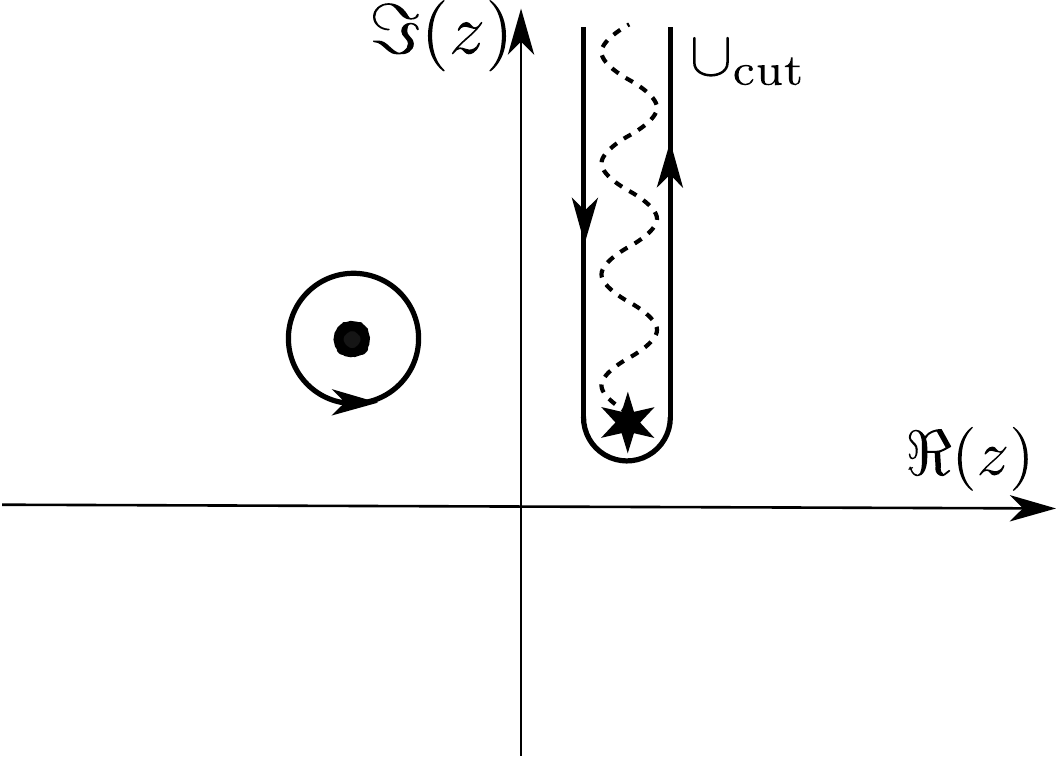}
  \caption{Singularity map of the integrand of (\ref{eq:bracketintegraldef}) in the $z$ UHP (left) and contour
  deformation around the pole and branch-cut (right).}
\label{fig:contour_trick}
\end{figure}

By definition of the Cauchy bracket, and the integral $I_{+ -}^{\varepsilon}$,
we have
\begin{eqnarray}
  \left[ \frac{I_{+ -}^{\varepsilon}}{K_{- \circ}} \right]_{+ \circ}
  \!\!\!\!(\tmmathbf{\alpha}) & = & \int_{\mathcal{A}^b} \Psi (z, \alpha_1, \alpha_2)
  I_{+ -}^{\varepsilon} (z, \alpha_2) \mathd z,  \label{eq:bracketintegraldef}
\end{eqnarray}
where
\begin{eqnarray*}
  \Psi (z, \alpha_1, \alpha_2) & = & \frac{1}{2 \pi i (z - \alpha_1) K_{- \circ}
  (z, \alpha_2)} \cdot
\end{eqnarray*}
It is interesting to note that the singularities of the integrand of
(\ref{eq:bracketintegraldef}) in the $z$ UHP are exclusively those of $\Psi (z,
\alpha_1, \alpha_2)$ since \ $I_{+ -}^{\varepsilon} (z, \alpha_2)$ is analytic
there. Hence, in the $z$ UHP, the integrand of (\ref{eq:bracketintegraldef}) has one simple pole at $z =
\alpha_1$ and one branch point at $z = \sqrt{k^2 - \alpha_2^2}$, with a branch-cut going vertically upwards as depicted in Figure~\ref{fig:contour_trick}
(left).

It is hence possible to deform the contour from $\mathcal{A}^b$ to a contour
$\cup_{\tmop{cut}}$ surrounding the branch cut. Assuming\footnote{Since we are
interested in the behaviour of this bracket for fixed $\alpha_1$ as $|
\alpha_2 | \rightarrow \infty$, we can make this assumption without loss of
generality.} for now that $\alpha_1 \neq \sqrt{k^2-\alpha_2^2}$, the pole at $z
= \alpha_1$ is picked up in the process and its contribution must be accounted
for, see Figure~\ref{fig:contour_trick} (right).

Noting that
\begin{eqnarray*}
  2 \pi i \underset{z = \alpha_1}{\tmop{Res}} (\Psi (z, \alpha_1, \alpha_2) I_{+
  -}^{\varepsilon} (z, \alpha_2)) & = & \frac{I_{+ -}^{\varepsilon} (\alpha_1,
  \alpha_2)}{K_{- \circ} (\alpha_1, \alpha_2)},
\end{eqnarray*}
we can write
\begin{eqnarray}
  \left[ \frac{I_{+ -}^{\varepsilon}}{K_{- \circ}} \right]_{+ \circ}
  \!\!\!\!(\tmmathbf{\alpha}) & = & \frac{I_{+ -}^{\varepsilon} (\alpha_1,
  \alpha_2)}{K_{- \circ} (\alpha_1, \alpha_2)} + I_{\tmop{cut}} (\alpha_1,
  \alpha_2),  \label{eq:bracket decomposition}
\end{eqnarray}
where
\begin{eqnarray*}
  I_{\tmop{cut}} (\alpha_1, \alpha_2) & = & \int_{\cup_{\tmop{cut}}} \Psi (z,
  \alpha_1, \alpha_2) I_{+ -}^{\varepsilon} (z, \alpha_2) \mathd z .
\end{eqnarray*}
Now using the fact that $\Psi$ changes sign across the cut, and that it is equal
to zero at the branch point, we can rewrite this integral in the slightly
simpler form
\begin{eqnarray*}
  I_{\tmop{cut}} (\alpha_1, \alpha_2) & = & 2 i \int_0^{\infty} \Psi \left(
  \sqrt{k^2 - \alpha_2^2} + it, \alpha_1, \alpha_2 \right) I_{+
  -}^{\varepsilon} \left( \sqrt{k^2 - \alpha_2^2} + it, \alpha_2 \right)
  \mathd t ,
\end{eqnarray*}
where $\Psi$ is only evaluated on the right side of its cut. 

\begin{figure}[h]
\centering
\includegraphics[width=0.5\textwidth]{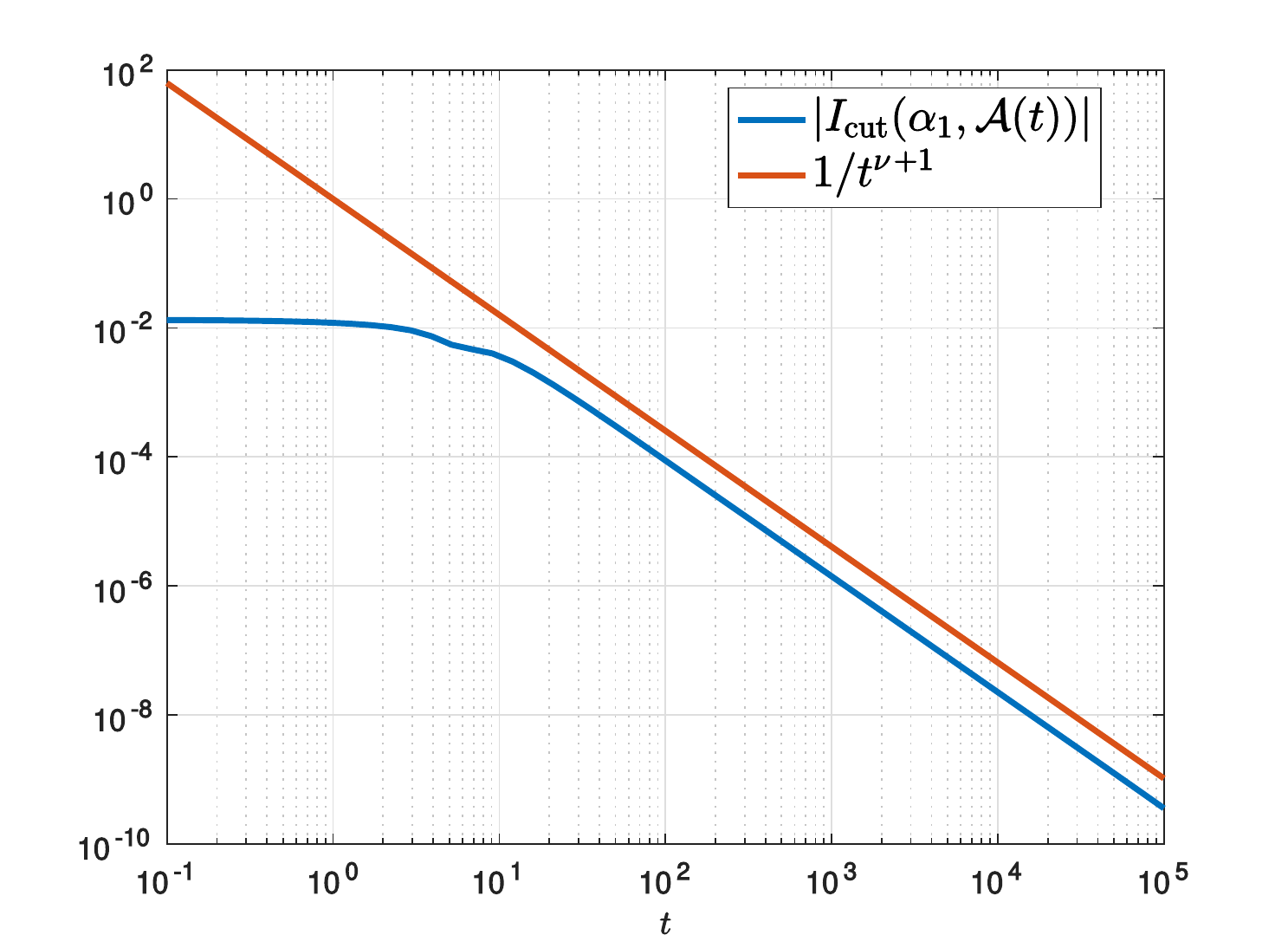}
  \caption{Log log plot of the absolute value of $I_{\tmop{cut}}$ for
  $\alpha_1 = 1 + 3 i$, and $\alpha_2$ tends to infinity on the contour
  $\mathcal{A}$. Here, $\mathcal{A}$ is parametrised by a real parameter $t$,
  hence the $\mathcal{A} (t)$ notation. $\mathcal{A} (t)$ here should be
  understood to take complex values, but be such that $\mathcal{A} (t) \sim t$
  as $t \rightarrow \infty$. The decay is compared to that of $1 / t^{\RED{\nu} +
  1}$. }
\label{fig:decayIcut}
\end{figure}

At this stage, it is useful to note that by (\ref{eq:defKminuso}), we have
\begin{eqnarray}
  K_{- \circ} (\alpha_1, \alpha_2) & \underset{| \alpha_2 | \rightarrow
  \infty}{\overset{\alpha_1 \tmop{fixed}}{=}} & \mathcal{O} (| \alpha_2 |^{- 1
  / 2}) .  \label{eq:Kminusobehaviour}
\end{eqnarray}
Using this and the asymptotic result
(\ref{eq:asymp-behav-3}), we find that the first term in the RHS of
(\ref{eq:bracket decomposition}) behaves like
\begin{eqnarray}
  \frac{I_{+ -}^{\varepsilon} (\alpha_1, \alpha_2)}{K_{- \circ} (\alpha_1,
  \alpha_2)} & \underset{| \alpha_2 | \rightarrow \infty}{\overset{\alpha_1
  \tmop{fixed}}{=}} & \mathcal{O} \left( \frac{1}{| \alpha_2 |} \right), \label{eqraph:behavpole}
\end{eqnarray}
which satisfies the condition (\ref{eq:condition}). Moreover, we
show in Appendix \ref{App:Icut} that we also have
\begin{eqnarray}
   I_{\tmop{cut}} (\alpha_1, \alpha_2)  & \underset{| \alpha_2 | \rightarrow
  \infty}{\overset{\alpha_1 \tmop{fixed}}{=}} & \mathcal{O} \left( \frac{1}{|
  \alpha_2 |^{\RED{\nu} + 1}} \right), \label{eqraph:behavIcut}
\end{eqnarray}
which also satisfies the condition (\ref{eq:condition}). Numerical evaluation
of this integral confirms this finding, as illustrated on Figure
\ref{fig:decayIcut}.

In conclusion, using (\ref{eq:bracket decomposition}), (\ref{eqraph:behavpole}) and (\ref{eqraph:behavIcut}), it is quite clear that as $| \alpha_2 | \rightarrow \infty$ on
$\mathcal{A}_2$, we have
\begin{eqnarray*}
  \left[ \frac{I_{+ -}^{\varepsilon}}{K_{- \circ}} \right]_{+ \circ}
  \!\!\!\!(\alpha_1, \alpha_2) & = & \mathcal{O} \left( \frac{1}{| \alpha_2 |}
  \right),
\end{eqnarray*}
meaning that the condition (\ref{eq:condition}) is fulfilled for $\gamma = 1$.
We can hence conclude this section by saying that
\begin{eqnarray*}
  \mathcal{I} [I_{+ -}^{\varepsilon}] & \underset{| \alpha_2 | \rightarrow
  \infty}{\longrightarrow} & 0.
\end{eqnarray*}
\section{Significance and perspectives}\label{sec:conclusion}


\RED{In this work we have introduced an explicit canonical integral, $I_{+-}^{\varepsilon}(\alpha_1,\alpha_2)$, that is very similar to the key unknown function $G_{+-}(\alpha_1,\alpha_2)$ in the quarter-plane problem, in the sense that they share the same asymptotic behaviour as $|\alpha_{1,2}|\to \infty$, and they share the same crucial property $\mathcal{I} [G_{+ -}] \rightarrow 0$ and $\mathcal{I} [I_{+ -}^{\varepsilon}] \to 0$. The latter was used in \cite{Assier2019} to prove the two
formulae (\ref{eq:radlowoccurence}) (including Radlow's ansatz) and (\ref{eq:compatibilityeq}) (coined the
\tmem{compatibility equation}).}

The exact form of the integral $I_{+ -}^{\varepsilon}$ has the
potential to be of great assistance to the design of a scheme to accurately
approximate the key unknown function $G_{+ -}$. One could for example consider
an approximation of the type
\begin{eqnarray*}
  G^{(N)}_{+ -} (\tmmathbf{\alpha}) & = & \mathcal{C}I_{+ -}^{\varepsilon}
  (\tmmathbf{\alpha}) T_{+ -} (\tmmathbf{\alpha}) \left( 1 + \sum_{j = 1}^N
  g_{+ -}^{(j)} (\tmmathbf{\alpha}) \right),
\end{eqnarray*}
where $\mathcal{C}$ is a constant, $T_{+ -} (\tmmathbf{\alpha})$ is a bounded function
(but not decaying to zero) at infinity, and the functions $g_{+ -}^{(j)}
(\tmmathbf{\alpha})$ are a set of simple functions (possibly \RED{with} simple poles at
given locations, but with unknown residues) that decay to zero at infinity. We
must choose $T_{+ -}$ and $g_{+ -}^{(j)}$ to be analytic in
$\mathcal{D}_{+ -}$. Note that the aim here is similar to that of \cite{Abrahams2000} say, \RED{where} functions analytic in a half-plane \RED{were approximated} via Pad\'e approximants.

We could for example choose the $T_{+ -}$ function to take the form
\begin{eqnarray*}
  T_{+ -} (\tmmathbf{\alpha}) & = & \frac{\mathcal{L} \alpha_1 -
  \alpha_2}{\alpha_1 - \alpha_2},
\end{eqnarray*}
for some unknown constant $\mathcal{L}$, while the $g_{+ -}^{(j)}$ could be
chosen to take the form
\begin{eqnarray*}
  g_{+ -}^{(j)} & = & \frac{\mathcal{R}^{(j)}}{(\alpha_1 - a_1^{(j)})
  (\alpha_2 - a_2^{(j)})},
\end{eqnarray*}
for some specified $(a_1^{(j)}, a_2^{(j)}) \in \mathcal{D}_{- +}$, and some
unknown residues $\mathcal{R}^{(j)}$.

\RED{The correcting function $T_{+-}$ and the constant $\mathcal{C}$ are included to compensate for the fact that the pre-factors (depending on $\varepsilon$) in the asymptotic results (\ref{eq:asymp-behav-1})--(\ref{eq:asymp-behav-3}) on $I_{+-}^{\varepsilon}$ may not be exactly equal to those of $G_{+-}$.}

For a given $N$, we will hence have $N + 2$ unknowns: $(\mathcal{C},
\mathcal{L}, \mathcal{R}^{(1)}, \ldots, \mathcal{R}^{(N)})$, which will be
determined by ensuring that the compatibility equation is satisfied at a set
of $N + 2$ collocation points. The implementation of such scheme is beyond the
scope of the present work and will constitute the basis of further investigations by the authors.

\paragraph{Acknowledgements} Both authors would like to thank the Isaac Newton Institute for Mathematical Sciences, Cambridge, for support and hospitality during the programme ``Bringing pure and applied analysis together via the Wiener--Hopf technique, its generalisations and applications'' where some work on this paper was undertaken. This work was supported by the EPSRC grant EP/R014604/1. Abrahams also acknowledges the support of UKRI/EPSRC grant EP/K032208/.

{
\clearpage
\bibliographystyle{plain}
\bibliography{biblio}
}

\appendix
\section{Proof of the asymptotic behaviour of $I_{+
-}^{\varepsilon}$}\label{app:app1}

\subsection{Proof of the asymptotic form (\ref{eq:asymp-behav-1})}

This is the simplest of the three different cases to be considered. We will
consider that both $| \alpha_1 |$ and $| \alpha_2 | \rightarrow \infty$ within
$\mathcal{D}_{+ -}$, in such a way that there exists an $m > 0$ such that $|
\alpha_2 | = m | \alpha_1 |$. As discussed in Section~\ref{sec:asymptoticbehaviour}, since both $\alpha_1$ and $\alpha_2$ are large, we can write $\alpha_{1,2}=|\alpha_{1,2}|e^{i\phi_{1,2}}$ for $\phi_1\in(0,\pi)$ and $\phi_2\in(-\pi,0)$. In this case, \RED{starting from (\ref{eq:initialdefIrevision})}, we have
\begin{eqnarray*}
  \PUR{I_{+ -}^{\varepsilon} (\alpha_1, \alpha_2)} 
  & \overset{\alpha_1 = | \alpha_1 | e^{i \phi_1}}{\underset{\alpha_2 = |
  \alpha_2 | e^{i \phi_2}}{=}} & \int_{\frac{3 \pi}{2}}^{2 \pi} \frac{\RED{f}
  \left( \frac{\pi}{2}, \varphi \right)}{(- i (| \alpha_1 | e^{i \phi_1} \cos
  (\varphi) + | \alpha_2 | e^{i \phi_2} \sin (\varphi) + i
  \varepsilon))^{\RED{\nu} + 3 / 2}} \mathd \varphi\\
  & \underset{| \alpha_2 | = m | \alpha_1 |}{=} & \frac{1}{| \alpha_1
  |^{\RED{\nu} + 3 / 2}} \int_{\frac{3 \pi}{2}}^{2 \pi} \frac{\RED{f} \left(
  \frac{\pi}{2}, \varphi \right)}{\left( - i \left( e^{i \phi_1} \cos
  (\varphi) + m e^{i \phi_2} \sin (\varphi) + \frac{i \varepsilon}{| \alpha_1
  |} \right) \right)^{\RED{\nu} + 3 / 2}} \mathd \varphi\\
  & \underset{| \alpha_1 | \rightarrow \infty}{\sim} & \frac{1}{| \alpha_1
  |^{\RED{\nu} + 3 / 2}} \int_{\frac{3 \pi}{2}}^{2 \pi} \frac{\RED{f} \left(
  \frac{\pi}{2}, \varphi \right)}{(- i (e^{i \phi_1} \cos (\varphi) + m e^{i
  \phi_2} \sin (\varphi)))^{\RED{\nu} + 3 / 2}} \mathd \varphi\\
  & \underset{| \alpha_1 | \rightarrow \infty}{\sim} & \frac{1}{| \alpha_1
  |^{\RED{\nu} + 3 / 2}} I_{+ -}^0 (e^{i \phi_1}, me^{i \phi_2})\underset{| \alpha_1 | \rightarrow \infty}{=} \mathcal{O} \left(
  \frac{1}{| \alpha_1 |^{\RED{\nu} + 3 / 2}} \right),
\end{eqnarray*}
since the quantity $e^{i \phi_1} \cos (\varphi) + m e^{i \phi_2} \sin
(\varphi)$ can never be equal to zero and hence the last integral is well
defined and independent of $| \alpha_1 |$.

This last statement can be proven as follows. Let us assume that $e^{i \phi_1}
\cos (\varphi) + m e^{i \phi_2} \sin (\varphi)$=0 for some $\varphi \in \left[
\frac{3 \pi}{2}, 2 \pi \right]$. If $\varphi = 3 \pi / 2$, then we have $- m
e^{i \phi_2} = 0$, which is impossible since $m > 0$. If $\varphi = 2 \pi$,
then we have $e^{i \phi_1} = 0,$ which is impossible. Hence our quantity
cannot be zero at the end points of the integration domain. We can hence assume
that $\cos (\varphi) \neq 0$ and $\sin (\varphi) \neq 0$ and rewrite the
equality as $e^{i (\phi_1 - \phi_2)} = - m \tan (\varphi) > 0$. Now taking the
imaginary part on both sides, we get $\sin (\phi_1 - \phi_2) = 0$, implying
that $\phi_1 = \phi_2 + n \pi$ for some $n \in \mathbb{Z}$. Clearly, from the
restriction on $\phi_{1, 2}$, we can only have $n = - 1, 0, 1, 2$. For $n = -
1, 1$, we would get $- 1 = - m \tan (\varphi) > 0$, which is impossible. Hence
we have $n = 0$ or $2$, i.e.\ $\phi_1 = \phi_2$ or $\phi_1 = \phi_2 + 2
\pi$. Because of the restriction on $\phi_{1, 2}$ this imposes $\phi_{1, 2} =
0$ or $\phi_1 = \pi$ and $\phi_2 = - \pi$, but these values are excluded
according to the restriction on $\phi_{1, 2}$\RED{, which contradicts our initial assumption and proves that $e^{i \phi_1} \cos (\varphi) + m e^{i \phi_2} \sin
(\varphi)\neq 0$}.

\subsection{Proof of the asymptotic form (\ref{eq:asymp-behav-2})}\label{sec:I1I2}

We have
\begin{eqnarray*}
  I_{+ -}^{\varepsilon} (\alpha_1, \alpha_2) & = & \int_{\frac{3 \pi}{2}}^{2
  \pi} \frac{\RED{f} \left( \frac{\pi}{2}, \varphi \right)}{(- i (\alpha_1 \cos
  (\varphi) + \alpha_2 \sin (\varphi) + i \varepsilon))^{\RED{\nu} + 3 / 2}}
  \mathd \varphi .
\end{eqnarray*}
We can see that as $| \alpha_1 | \rightarrow \infty$, the denominator is
dominated by the term involving $\alpha_1 \cos (\varphi)$ for all $\varphi$,
except when $\varphi \approx \frac{3 \pi}{2}$, where $\cos (\varphi)$
approaches zero. We can hence split the integral into two parts as $I_{+
-}^{\varepsilon} (\alpha_1, \alpha_2) = I_1^{\varepsilon, \delta} (\alpha_1,
\alpha_2) + I_2^{\varepsilon, \delta} (\alpha_1, \alpha_2)$, where
\begin{eqnarray}
  I_1^{\varepsilon, \delta} (\alpha_1, \alpha_2) & = & \int_{\frac{3
  \pi}{2}}^{\frac{3 \pi}{2} + \delta} \frac{\RED{f} \left( \frac{\pi}{2}, \varphi
  \right)}{(- i (\alpha_1 \cos (\varphi) + \alpha_2 \sin (\varphi) + i
  \varepsilon))^{\RED{\nu} + 3 / 2}} \mathd \varphi, \\
  I_2^{\varepsilon, \delta} (\alpha_1, \alpha_2) & = & \int_{\frac{3 \pi}{2} +
  \delta}^{2 \pi} \frac{\RED{f} \left( \frac{\pi}{2}, \varphi \right)}{(- i
  (\alpha_1 \cos (\varphi) + \alpha_2 \sin (\varphi) + i \varepsilon))^{\RED{\nu}
  + 3 / 2}} \mathd \varphi, 
\end{eqnarray}
for some small $\delta = \delta (| \alpha_1 |) > 0$, chosen such that $|
\alpha_1 | \delta \rightarrow \infty$ as $| \alpha_1 | \rightarrow \infty$ and
$\delta \rightarrow 0$ as $| \alpha_1 | \rightarrow \infty$. Let us select $\delta = 1 / \sqrt{| \alpha_1 |}$ as it will prove to work. Upon making the change of
variable $\psi = \varphi - \frac{3 \pi}{2}$, these integrals become
\begin{eqnarray*}
  I_1^{\varepsilon, \delta} (\alpha_1, \alpha_2) & = & \int_0^{\delta}
  \frac{\RED{f} \left( \frac{\pi}{2}, \frac{3 \pi}{2} + \psi \right)}{(- i
  (\alpha_1 \sin (\psi) - \alpha_2 \cos (\psi) + i \varepsilon))^{\RED{\nu} + 3 /
  2}} \mathd \psi,\\
  I_2^{\varepsilon, \delta} (\alpha_1, \alpha_2) & = &
  \int_{\delta}^{\frac{\pi}{2}} \frac{\RED{f} \left( \frac{\pi}{2}, \psi + \frac{3
  \pi}{2} \right)}{(- i (\alpha_1 \sin (\psi) - \alpha_2 \cos (\psi) + i
  \varepsilon))^{\RED{\nu} + 3 / 2}} \mathd \psi .
\end{eqnarray*}
For the first integral, we have
\begin{eqnarray}
  I_1^{\varepsilon, \delta} (\alpha_1, \alpha_2) & \underset{\delta
  \rightarrow 0}{\sim} & \int_0^{\delta} \frac{\RED{f} \left( \frac{\pi}{2},
  \frac{3 \pi}{2} \right)}{(- i (\alpha_1 \psi - \alpha_2 + i
  \varepsilon))^{\RED{\nu} + 3 / 2}} \mathd \psi .  \label{eq:plvsdvsr}
\end{eqnarray}
Since $| \alpha_1 | \rightarrow \infty$, remember that we can write $\alpha_1
= | \alpha_1 | e^{i \phi_1}$ with $\phi_1\in(0,\pi)$, and making the change of variable $\theta = |
\alpha_1 | \psi$ in (\ref{eq:plvsdvsr}) we get
\begin{eqnarray*}
  I_1^{\varepsilon, \delta} (\alpha_1, \alpha_2) & \underset{}{\sim} &
  \int_0^{| \alpha_1 | \delta} \frac{\RED{f} \left( \frac{\pi}{2}, \frac{3 \pi}{2}
  \right)}{(- i (\theta e^{i \phi_1} - \alpha_2 + i \varepsilon))^{\RED{\nu} + 3 /
  2}}  \frac{\mathd \theta}{| \alpha_1 |}\\
  & \underset{\delta | \alpha_1 | \rightarrow \infty}{\sim} & \frac{\RED{f}
  \left( \frac{\pi}{2}, \frac{3 \pi}{2} \right)}{| \alpha_1 |} 
  \int_0^{\infty} \frac{1}{(- i (\theta e^{i \phi_1} - \alpha_2 + i
  \varepsilon))^{\RED{\nu} + 3 / 2}} \mathd \theta .
\end{eqnarray*}
The latter integral can be recast in the form $\int_0^{\infty} \frac{\mathd
\theta}{(A - B \theta)^{\lambda}}$, for $A = i \alpha_2 + \varepsilon$, $B = i
e^{i \phi_1}$ and $\lambda = \RED{\nu} + 3 / 2$. Since we have $\tmop{Re} (A) >
0$, $A / B \nin \mathbb{R}$ and $\lambda > 1$, this integral can easily be shown to
be equal to $\frac{A^{1 - \lambda}}{B (1 - \lambda)}$, and hence
\begin{eqnarray}
  \PUR{I_1^{\varepsilon, \delta} (\alpha_1, \alpha_2)} & \underset{| \alpha_1 |
  \rightarrow \infty}{\overset{\alpha_2 \tmop{fixed}}{\sim}} & 
\frac{1}{\alpha_1} \times \frac{2 i \RED{f} \left( \frac{\pi}{2}, \frac{3
  \pi}{2} \right)}{(i \alpha_2 + \varepsilon)^{\RED{\nu} + 1 / 2} (2 \RED{\nu} + 1)} \cdot
  \label{eq:behaveI1raph}
\end{eqnarray}

For the second integral, we can use the fact that $| \alpha_1 \sin (\psi) | >
| \alpha_1 \sin (\delta) | \rightarrow \infty$ to simplify the denominator and
obtain
\begin{eqnarray*}
  \PUR{I_2^{\varepsilon, \delta} (\alpha_1, \alpha_2)} & 
\RED{\underset{| \alpha_1 |
  \rightarrow \infty}{\overset{\alpha_1 = | \alpha_1 | e^{i \phi_1}}{\sim}}} &
  \frac{1}{| \alpha_1
  |^{\RED{\nu} + 3 / 2}} \int_{\delta}^{\frac{\pi}{2}} \frac{\RED{f} \left(
  \frac{\pi}{2}, \psi + \frac{3 \pi}{2} \right)}{(- ie^{i \phi_1} \sin
  (\psi))^{\RED{\nu} + 3 / 2}} \mathd \psi .
\end{eqnarray*}
Hence, since by Lemma \ref{lem:eigenfunctionproperties}, $0 \leqslant \RED{f}
\left( \frac{\pi}{2}, \psi + \frac{3 \pi}{2} \right) \leqslant \RED{f} \left(
\frac{\pi}{2}, \frac{3 \pi}{2} \right)$, and because $| \sin (\delta) |
\leqslant | \sin (\psi) |$, we have
\begin{eqnarray*}
  | I_2^{\varepsilon, \delta} (\alpha_1, \alpha_2) | & \underset{| \alpha_1 |
  \rightarrow \infty}{\sim} & \frac{1}{| \alpha_1 |^{\RED{\nu} + 3 / 2}}
  \int_{\delta}^{\frac{\pi}{2}} \frac{\RED{f} \left( \frac{\pi}{2}, \psi + \frac{3
  \pi}{2} \right)}{| - ie^{i \phi_1} \sin (\psi) |^{\RED{\nu} + 3 / 2}} \mathd
  \psi\\
  & \leqslant & \frac{\frac{\pi}{2} \RED{f} \left( \frac{\pi}{2}, \frac{3 \pi}{2}
  \right)}{(| \sin (\delta) | | \alpha_1 |)^{\RED{\nu} + 3 / 2}} \
  \underset{\delta \rightarrow 0}{\sim} \ \frac{\frac{\pi}{2} \RED{f}
  \left( \frac{\pi}{2}, \frac{3 \pi}{2} \right)}{(\delta | \alpha_1 |)^{\RED{\nu}
  + 3 / 2}} \ \underset{\delta = \frac{1}{\sqrt{\alpha_1}}}{=}
  \mathcal{O} \left( \frac{1}{| \alpha_1 |^{\frac{\RED{\nu}}{2} + \frac{3}{4}}}
  \right) = \smallO \left( \frac{1}{| \alpha_1 |} \right)
\end{eqnarray*}
since $\frac{\RED{\nu}}{2} + \frac{3}{4} \approx 1.15 > 1$. Hence, overall,
$I_2^{\varepsilon, \delta}$ can be neglected to leading order and, using
(\ref{eq:behaveI1raph}), we obtain
\begin{eqnarray*}
  I_{+ -}^{\varepsilon} (\alpha_1, \alpha_2) \underset{| \alpha_1 |
  \rightarrow \infty}{\overset{\alpha_2 \tmop{fixed}}{\sim}} \frac{\Lambda_1
  (\alpha_2, \varepsilon)}{\alpha_1}, & \tmop{where} & \Lambda_1 (\alpha_2,
  \varepsilon) = \frac{2 i\RED{f} \left( \frac{\pi}{2}, \frac{3 \pi}{2}
  \right)}{(\varepsilon + i \alpha_2)^{\RED{\nu} + 1 / 2}  (1 + 2 \RED{\nu})},
\end{eqnarray*}
as required.

\subsection{Proof of the asymptotic form (\ref{eq:asymp-behav-3})}\label{sec:I3I4}

We have
\begin{eqnarray*}
  I_{+ -}^{\varepsilon} (\alpha_1, \alpha_2) & = & \int_{\frac{3 \pi}{2}}^{2
  \pi} \frac{\RED{f} \left( \frac{\pi}{2}, \varphi \right)}{(- i (\alpha_1 \cos
  (\varphi) + \alpha_2 \sin (\varphi) + i \varepsilon))^{\RED{\nu} + 3 / 2}}
  \mathd \varphi ,
\end{eqnarray*}
and we can see that as $| \alpha_2 | \rightarrow \infty$, the denominator is
dominated by the term involving $\alpha_2 \sin (\varphi)$ for all $\varphi$,
except when $\varphi \sim 2 \pi$, in which region $\sin (\varphi)$ approaches zero. We
should hence split the integral into two parts as $I_{+ -}^{\varepsilon}
(\alpha_1, \alpha_2) = I_3^{\varepsilon, \delta} (\alpha_1, \alpha_2) +
I_4^{\varepsilon, \delta} (\alpha_1, \alpha_2)$, where
\begin{eqnarray}
  I_3^{\varepsilon, \delta} (\alpha_1, \alpha_2) & = & \int_{2 \pi -
  \delta}^{2 \pi} \frac{\RED{f} \left( \frac{\pi}{2}, \varphi \right)}{(- i
  (\alpha_1 \cos (\varphi) + \alpha_2 \sin (\varphi) + i \varepsilon))^{\RED{\nu}
  + 3 / 2}} \mathd \varphi,  \label{eq:defI3}\\
  I_4^{\varepsilon, \delta} (\alpha_1, \alpha_2) & = & \int_{\frac{3
  \pi}{2}}^{2 \pi - \delta} \frac{\RED{f} \left( \frac{\pi}{2}, \varphi
  \right)}{(- i (\alpha_1 \cos (\varphi) + \alpha_2 \sin (\varphi) + i
  \varepsilon))^{\RED{\nu} + 3 / 2}} \mathd \varphi,  \label{eq:defI4}
\end{eqnarray}
for some small $\delta = \delta (| \alpha_2 |) > 0$, chosen such that $|
\alpha_2 | \delta \rightarrow \infty$ as $| \alpha_2 | \rightarrow \infty$ and
$\delta \rightarrow 0$ as $| \alpha_2 | \rightarrow \infty$. For specificity,
let us select $\delta = 1 / | \alpha_2 |^{1 / 4}$. Upon making the change of
variable $\psi = 2 \pi - \varphi$, these integrals become
\begin{eqnarray*}
  I_3^{\varepsilon, \delta} (\alpha_1, \alpha_2) & = & \int_0^{\delta}
  \frac{\RED{f} \left( \frac{\pi}{2}, 2 \pi - \psi \right)}{(- i (\alpha_1 \cos
  (\psi) - \alpha_2 \sin (\psi) + i \varepsilon))^{\RED{\nu} + 3 / 2}} \mathd
  \psi,\\
  I_4^{\varepsilon, \delta} (\alpha_1, \alpha_2) & = &
  \int_{\delta}^{\frac{\pi}{2}} \frac{\RED{f} \left( \frac{\pi}{2}, 2 \pi - \psi
  \right)}{(- i (\alpha_1 \cos (\psi) - \alpha_2 \sin (\psi) + i
  \varepsilon))^{\RED{\nu} + 3 / 2}} \mathd \psi .
\end{eqnarray*}
For the first integral, using Lemma \ref{lem:eigenfunctionproperties}, we have
\begin{eqnarray}
  I_3^{\varepsilon, \delta} (\alpha_1, \alpha_2) & \underset{\delta
  \rightarrow 0}{\sim} & \int_0^{\delta} \frac{\beta \sqrt{\psi}}{(- i
  (\alpha_1 - \alpha_2 \psi + i \varepsilon))^{\RED{\nu} + 3 / 2}} \mathd \psi 
  \label{eq:asfsadvoasni}\\
  & \underset{\theta = | \alpha_2 | \psi}{\overset{\alpha_2 = | \alpha_2 |
  e^{i \phi_2}}{=}} & \int_0^{| \alpha_2 | \delta} \frac{\beta
  \sqrt{\frac{\theta}{| \alpha_2 |}}}{(- i (\alpha_1 - e^{i \phi_2} \theta + i
  \varepsilon))^{\RED{\nu} + 3 / 2}} \frac{\mathd \theta}{| \alpha_2 |}
  \nonumber\\
  & \underset{| \alpha_2 | \delta \rightarrow \infty}{\sim} & \frac{\beta}{|
  \alpha_2 |^{3 / 2}} \int_0^{\infty} \frac{\sqrt{\theta}}{(- i (\alpha_1 -
  e^{i \phi_2} \theta + i \varepsilon))^{\RED{\nu} + 3 / 2}} \mathd \theta . \nonumber
\end{eqnarray}
The latter integral can be recast in the form $\int_0^{\infty}
\frac{\sqrt{\theta}}{(A - B \theta)^{\lambda}} \mathd \theta$, with $A = (- i
\alpha_1 + \varepsilon)$, $B = - i e^{i \phi_2}$ and $\lambda = \RED{\nu} + 3 /
2$. Since we have $\tmop{Re} (A) > 0$, $A / B \nin \mathbb{R}$ and $\lambda >
3 / 2$, this integral can be shown to be equal to $\frac{A^{- \lambda}
\sqrt{\pi} \Gamma (\lambda - 3 / 2)}{2 (- B / A)^{3 / 2} \Gamma (\lambda)}$
(see {\cite{GRADSHTEYN1980211}} p285, 3.194.3.), and hence
\begin{eqnarray}
  \PUR{I_3^{\varepsilon, \delta} (\alpha_1, \alpha_2)} & \underset{| \alpha_2 |
  \rightarrow \infty}{\overset{\alpha_1 \tmop{fixed}}{\sim}} & 
\frac{1}{\alpha_2^{3 / 2}} \times \frac{\beta (- i)^{3 / 2}
  \sqrt{\pi} \Gamma (\RED{\nu})}{2 (\varepsilon - i \alpha_1)^{\RED{\nu}} \Gamma
  (\RED{\nu} + 3 / 2)} \cdot  \label{eq:behaveI3epsdel}
\end{eqnarray}
For the second integral, since $| \alpha_2 \sin (\psi) | > | \alpha_2 \sin
(\delta) | \rightarrow \infty$, we have
\begin{eqnarray*}
  \PUR{I_4^{\varepsilon, \delta} (\alpha_1, \alpha_2)} & \underset{| \alpha_2 |
  \rightarrow \infty}{\sim} &
\frac{1}{(i \alpha_2)^{\RED{\nu} + 3 / 2}} \int_{\delta}^{\frac{\pi}{2}}
  \frac{\RED{f} \left( \frac{\pi}{2}, 2 \pi - \psi \right)}{(\sin (\psi))^{\RED{\nu} +
  3 / 2}} \mathd \psi .
\end{eqnarray*}
Hence, using Lemma \ref{lem:eigenfunctionproperties}, we have
\begin{eqnarray*}
  | I_4^{\varepsilon, \delta} (\alpha_1, \alpha_2) | & < & \frac{\RED{f} \left(
  \frac{\pi}{2}, \frac{3 \pi}{2} \right) \frac{\pi}{2}}{| \alpha_2 \sin
  (\delta) |^{\RED{\nu} + 3 / 2}} \  \underset{\delta = | \alpha_2 |^{-
  1 / 4}}{=} \ \frac{\pi \RED{f} \left( \frac{\pi}{2}, \frac{3
  \pi}{2} \right)}{2 | \alpha_2  |^{\frac{3 \RED{\nu}}{4} + \frac{9}{8}}}
  = \smallO \left( \frac{1}{| \alpha_2 |^{3 / 2}} \right),
\end{eqnarray*}
since $\frac{3 \RED{\nu}}{4} + \frac{9}{8} \approx 1.72 > \frac{3}{2}$.

Hence, overall, $I_4^{\varepsilon, \delta}$ can be neglected to leading order
and, using (\ref{eq:behaveI3epsdel}), we obtain
\begin{eqnarray*}
  I_{+ -}^{\varepsilon} (\alpha_1, \alpha_2) \underset{| \alpha_2 |
  \rightarrow \infty}{\overset{\alpha_1 \tmop{fixed}}{\sim}} \frac{\Lambda_2
  (\alpha_1, \varepsilon)}{\alpha_2^{3 / 2}}, & \tmop{where} & \Lambda_2
  (\alpha_1, \varepsilon) = \frac{\beta e^{-3i\frac{\pi}{4}} \sqrt{\pi} \Gamma (\RED{\nu})}{2 (\varepsilon - i \alpha_1)^{\RED{\nu}} \Gamma (\RED{\nu} + 3 /
  2)},
\end{eqnarray*}
as required.

\section{Asymptotic behaviour of $I_{\tmop{cut}}$ }\label{App:Icut}

Since we are only interested in the behaviour of $I_{\tmop{cut}}$ when $|
\alpha_2 | \rightarrow \infty$ and $\alpha_2 \in \mathcal{A}_2$, we can
consider $\alpha_2$ to be real here. Let us assume that $\alpha_2 > 0$ and
$\alpha_2 \rightarrow \infty$; the $\alpha_2$ negative case can be dealt with in a similar fashion. Using the definitions of the functions $K_{-\circ}(\alpha_1,\alpha_2)$ and $\sqrt{k^2-\alpha_2^2}$ given in  \cite{Assier2019}, one can show that, on the right \RED{side} of the cut, $1 / K_{- \circ} \left( \sqrt{k^2 - \alpha_2^2} + it, \alpha_2
\right) = \sqrt{t} e^{\frac{3 i \pi}{4}}$ and that, as $\alpha_2 \rightarrow
\infty$, $\sqrt{k^2 - \alpha_2^2} \sim i \alpha_2$, we obtain
\begin{eqnarray*}
  I_{\tmop{cut}} (\alpha_1, \alpha_2) & = & 2 i \int_0^{\infty} \Psi \left(
  \sqrt{k^2 - \alpha_2^2} + it, \alpha_1, \alpha_2 \right) I_{+
  -}^{\varepsilon} \left( \sqrt{k^2 - \alpha_2^2} + it, \alpha_2 \right)
  \mathd t\\
  & \underset{\alpha_2 \rightarrow \infty}{\approx} & \frac{e^{\frac{3 i
  \pi}{4}}}{\pi} \int_0^{\infty} \frac{\sqrt{t}}{(i \alpha_2 + i t -
  \alpha_1)} \int_{\frac{3 \pi}{2}}^{2 \pi} \frac{\RED{f} \left( \frac{\pi}{2},
  \varphi \right) \mathd \varphi \mathd t}{(\alpha_2 e^{- i \varphi} + t \cos
  (\varphi) + \varepsilon)^{\RED{\nu} + 3 / 2}} \cdot
\end{eqnarray*}
Now make the substitution $t = \alpha_2 u$ to get
\begin{eqnarray*}
  \PUR{I_{\tmop{cut}} (\alpha_1, \alpha_2)} 
& \underset{\alpha_2 \rightarrow
  \infty}{\approx} &
\frac{e^{\frac{3 i \pi}{4}}}{\pi}  \frac{1}{\alpha_2^{\RED{\nu} + 1}}
  \int_0^{\infty} \frac{\sqrt{u}}{\left( i (u + 1) - \frac{\alpha_1}{\alpha_2}
  \right)} \int_{\frac{3 \pi}{2}}^{2 \pi} \frac{\RED{f} \left( \frac{\pi}{2},
  \varphi \right) \mathd \varphi \mathd u}{\left( e^{- i \varphi} + u \cos
  (\varphi) + \frac{\varepsilon}{\alpha_2} \right)^{\RED{\nu} + 3 / 2}} \cdot
\end{eqnarray*}
Now the denominators are never zero, even when neglecting the small terms
involving $1 / \alpha_2$, so we have
\begin{eqnarray}
  I_{\tmop{cut}} (\alpha_1, \alpha_2) & \underset{\alpha_2 \rightarrow
  \infty}{\approx} & \frac{e^{\frac{3 i \pi}{4}}}{\pi} 
  \frac{1}{\alpha_2^{\RED{\nu} + 1}} \int_0^{\infty} \frac{\sqrt{u}}{(i (u + 1))}
  \underbrace{\int_{\frac{3 \pi}{2}}^{2 \pi} \frac{\RED{f} \left( \frac{\pi}{2},
  \varphi \right) \mathd \varphi}{(e^{- i \varphi} + u \cos (\varphi))^{\RED{\nu}
  + 3 / 2}}}_{J (u)} \mathd u . \label{eq:convergentdoubleintegral}
\end{eqnarray}
We need to check that the double integral actually exists.
There are no convergence problems as $u \to 0$ and no singularity of the
integrand for $u \in \mathbb{R^+}$. So, we just need to ensure that the integral
converges at $\infty$. For this purpose, it is enough to study the behaviour
of $J (u)$ as $u \rightarrow \infty$.

The only possible issue occurs when $\cos (\varphi) \approx 0$; to resolve this, consider a small fixed constant $\delta > 0$, and write
\begin{eqnarray*}
  J (u) & = & \int_{\frac{3 \pi}{2}}^{\frac{3 \pi}{2} + \delta} \frac{\RED{f}
  \left( \frac{\pi}{2}, \varphi \right) \mathd \varphi}{(e^{- i \varphi} + u
  \cos (\varphi))^{\RED{\nu} + 3 / 2}} + \int_{\frac{3 \pi}{2} + \delta}^{2 \pi}
  \frac{\RED{f} \left( \frac{\pi}{2}, \varphi \right) \mathd \varphi}{(e^{- i
  \varphi} + u \cos (\varphi))^{\RED{\nu} + 3 / 2}}\\
  & = & J_1 (u) + J_2 (u)
\end{eqnarray*}
Let us start by studying $J_2 (u)$:
\begin{eqnarray*}
  | J_2 (u) | & = & \left| \int_{\frac{3 \pi}{2} + \delta}^{2 \pi} \frac{\RED{f}
  \left( \frac{\pi}{2}, \varphi \right) \mathd \varphi}{(e^{- i \varphi} + u
  \cos (\varphi))^{\RED{\nu} + 3 / 2}} \right|\\
  & \underset{u \rightarrow \infty}{\approx} & \left| \int_{\frac{3 \pi}{2} +
  \delta}^{2 \pi} \frac{\RED{f} \left( \frac{\pi}{2}, \varphi \right) \mathd
  \varphi}{(u \cos (\varphi))^{\RED{\nu} + 3 / 2}} \right|\\
  & < & \frac{1}{u^{\RED{\nu} + 3 / 2}}  \frac{\RED{f} \left( \frac{\pi}{2}, \frac{3
  \pi}{2} \right) \frac{\pi}{2}}{\delta^{\RED{\nu} + 3 / 2}} =\mathcal{O} \left(
  \frac{1}{u^{\RED{\nu} + 3 / 2}} \right) .
\end{eqnarray*}
Now consider the slightly more subtle case of $J_1$, and using Lemma
\ref{lem:eigenfunctionproperties}, we obtain
\begin{eqnarray*}
  \PUR{J_1 (u)} & = & \int_{\frac{3 \pi}{2}}^{\frac{3 \pi}{2} + \delta} \frac{\RED{f}
  \left( \frac{\pi}{2}, \varphi \right) \mathd \varphi}{(e^{- i \varphi} + u
  \cos (\varphi))^{\RED{\nu} + 3 / 2}}\\
  & \overset{z = \varphi - \frac{3 \pi}{2}}{\underset{\delta \ll 1}{\approx}}
  & \RED{f} \left( \frac{\pi}{2}, \frac{3 \pi}{2} \right) \int_0^{\delta}
  \frac{\mathd z}{(i + uz)^{\RED{\nu} + 3 / 2}}\\
  & \RED{\overset{x = u z}{\underset{\delta u \rightarrow \infty}{\approx}}} & \frac{\RED{f} \left(
  \frac{\pi}{2}, \frac{3 \pi}{2} \right)}{u} \int_0^{\infty} \frac{\mathd
  x}{(i + x)^{\RED{\nu} + 3 / 2}} =\mathcal{O} \left( \frac{1}{u} \right).
\end{eqnarray*}
Hence, since $\RED{\nu} + 3 / 2 > 1$, $J_2 (u)$ can be neglected to leading order,
and, overall $J (u) =\mathcal{O} (1 / u)$. This ensures that the double
integral in (\ref{eq:convergentdoubleintegral}) converges. Since it is
independent from $\alpha_2$, it is quite clear that
 $ I_{\tmop{cut}} = \mathcal{O} \left( \frac{1}{|\alpha_2|^{\RED{\nu} + 1}}
  \right)$,
as claimed in Section \ref{sec:proofstrategy}, and as confirmed numerically.

\end{document}